\documentclass[11pt]{amsart}

\usepackage{psfrag} 
\usepackage{epsfig}
\usepackage{amsmath}
\usepackage{graphicx}
\usepackage{amssymb}
\usepackage{amsthm}
\usepackage{graphics}
\usepackage{latexsym}

\newtheorem{theorem}{Theorem}

\newtheorem{lemma}{Lemma}
\newtheorem{remark}{Remark}

\newtheorem{corollary}{Corollary}
\newtheorem{proposition}{Proposition}
\newtheorem{fact}{Fact}
\newcommand{\al}{\alpha}
\newcommand{\ep}{\epsilon}
\newcommand{\pt}{\partial_t}
\newcommand{\pa}{\partial_\al}
\newcommand{\pbm}{\partial_{\beta_-}}
\newcommand{\pb}{\partial_{\beta}}
\newcommand{\pbp}{\partial_{\beta_+}}
\newcommand{\pT}{\partial_T}
\newcommand{\ptau}{\partial_\tau}
\newcommand{\bp}{\beta_+}
\newcommand{\bm}{\beta_-}
\newcommand{\bpm}{\beta_\pm}
\newcommand{\pz}{\Psi^z}

\newcommand{\py}{\Psi^y}

\newcommand{\pu}{\Psi^u}
\newcommand{\pd}{\Psi^d}
\newcommand{\putt}{{\Psi}^u}
\newcommand{\paa}{\Psi^a}
\newcommand{\patt}{{\Psi}^a}
\newcommand{\sob}{\mathfrak{H}}
\newcommand{\ab}{\bar{\al}}
\newcommand{\res}{\textrm{Res}}
\newcommand{\resp}{\res[\Psi]}

\newcommand{\third}{\frac{1}{3}}
\newcommand{\Kd}{\ep\pb+ \third \ep^3 \pb^3}
\newcommand{\Kt}{\ep \pb+\third \ep^3  \pb^3+\frac{2}{15} \ep^5 \pb^5}

\newcommand{\Kep}{K_{0,\ep}}
\newcommand{\Lep}{L_\ep}

\begin{document}

\title{Corrections to the KdV approximation for water waves}
\author{J. Douglas Wright}
\begin{abstract}
In order to investigate corrections to the common KdV approximation for
surface water waves in a canal, we derive modulation equations for the evolution of long
wavelength initial data.  We work in Lagrangian coordinates.
The equations which govern corrections to the KdV approximation
consist of linearized and
inhomogeneous KdV equations plus an inhomogeneous wave equation.  
These equations are explicitly solvable and   
we prove estimates showing that they do indeed give a
significantly better approximation than the KdV equation alone. 
\end{abstract}

\maketitle


\section{Introduction}

\begin{figure}
\begin{center}
  \psfragscanon
  \includegraphics{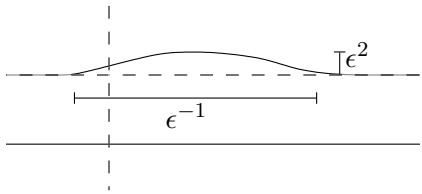}
  \caption{The long wave, small amplitude scaling.}
\label{scaling picture}
\end{center}
\end{figure}

It is often easier to write down a partial differential equation
which models a physical phenomena than it is to study solutions
of such an equation.  
Equations which model the evolution
of the surface of a fluid in a canal have been known since at least
the 19th century, however it has only been in recent years that questions 
of existence and uniqueness for general initial data have been answered
(see \cite{wu:97}, \cite{wu:99} and \cite{lannes.waterwaves:03}).  
Moreover, numerical simulations of water waves 
and similarly complex phenomena are frequently 
time consuming and challenging to implement.
Consequently, it can be quite difficult
to say much about the behavior of a general solution. And so
scientists often restrict their attention to limiting
cases---for instance, one may assume that solutions
are of long wavelength and small amplitude
(see Figure \ref{scaling picture}). 
Under such a supposition, a {\it modulation equation} 
may be (formally) derived.  In particular, one hopes that the modulation
equation:
\begin{itemize}
\item is well-posed,
\item is either explicitly solvable or easy to solve numerically and
\item captures the essential behavior of the original system.
\end{itemize}
Remarkably, many seemingly disparate
physical phenomena possess modulation equations of the same form.
For solutions of long wavelength,
Korteweg-de Vries (KdV) equations are often used
as modulation equations for a wide variety
of non-linear dispersive systems, including
the water wave equation, the Euler-Poisson
equations for plasma dynamics and the Fermi-Pasta-Ulam equation
for the interaction of particles in an infinite lattice.  

Despite the fact that modulation equations have been in use for over
a hundred years---the KdV equation was first proposed as a model for
water waves by Boussinesq in $1872$ and also by Korteweg and 
de Vries in $1895$---only recently have attempts been made to 
rigorously connect the behavior of the modulation equations to the
original physical problem.  In particular, 
through the work of Craig \cite{craig:85}, 
Kano and Nishida \cite{kano.etal:79} \cite{kano.etal:86}, 
Kalyakin \cite{kalyakin:89}, 
Schneider \cite{schneider:98}, Ben Youssef and Colin \cite{youssef.etal:00} 
and Schneider and Wayne \cite{schneider.etal:00},
\cite{schneider.etal:02}, the validity of 
KdV equations as a leading order approximation
to the evolution of long wavelength water waves and to a 
number of other dispersive partial differential 
equations has been established.

In many respects, the KdV equation is an ideal modulation equation;  
it is simple in form and explicitly solvable {\it via} the inverse
scattering transform.  Nevertheless, both experimentally and numerically
one observes deviations from the 
predictions of the KdV approximation.  In this paper
we derive a hierarchy of modulation equations which govern corrections
to the KdV model and also prove rigorously that
these higher order equations do indeed improve the accuracy of
the approximation.   While the correction is valid in general long
wavelength/small amplitude settings, heuristically the model is
set up to better approximate interactions between solitary
waves---both 
counterpropagating collisions and unidirectional interactions.  

Note that in \cite{wayne.etal:02}, as a case study, Wayne and Wright examined
higher order corrections to the KdV approximation to a Boussinesq
equation.  As the KdV equation is in some sense a universal approximation
for long waves, we expect that the equations for corrections
to this approximation will also be universal.
Indeed, our results 
show that the higher order corrections for the 
water wave equation are nearly identical to those 
for the Boussinesq equation.  Also, since a significant part of 
the work for the Boussinesq problem consists of showing that the
modulation equations have well-behaved solutions over the time
scales of interest, this is of use in tackling the water wave
problem.

\begin{figure}
\begin{center}
  \psfragscanon
  \includegraphics{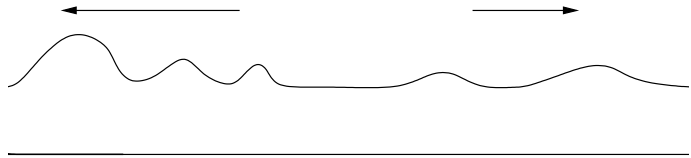} 
  \caption{Sketch of the right and left moving wave trains.}
\label{main result cartoon}
\end{center}
\end{figure}

We now describe our results in some detail.  
The equations of motion for a water wave in an infinitely long canal (commonly called
{\it the water wave equation}) are
\begin{equation}
\label{water wave}
\begin{split}
x_{tt}(1+x_\alpha) &+ y_\alpha(1+y_{tt}) = 0\\
y_t &= K(x,y) x_t\\
\alpha \in \mathbb{R},\, 
t\geq 0,\,& 
\left(x(\alpha,t),y(\alpha,t)\right) \in \mathbb{R}^2, 
\end{split}\tag{WW}
\end{equation}
where $K(x,y)$ is a complicated operator (see Section \ref{operator K 1}
for the definition of $K$) and $(\al + x(\al,t),y(\al,t))$ parameterizes
the free surface.  
According to the KdV approximation results of 
\cite{schneider.etal:00}, 
to the order of the approximation
long wavelength solutions of (\ref{water wave}) split up into two pieces, one
a right moving wave train and one a left moving wave train.  Each
of these wave trains evolves according to a KdV equation, and {\it there is
no interaction between the left and right moving pieces}.  That is,
for $0<\epsilon \ll 1$, if we scale amplitudes to be 
$O(\epsilon^2)$ ({\it i.e.} small) and wavelengths to be
$O(\epsilon^{-1})$ ({\it i.e.} long), then for 
times of $O(\epsilon^{-3})$, 
solutions to (\ref{water wave}) satisfy
\begin{equation}
\label{kdv approx for water wave}
-x_\alpha(\alpha,t)=\epsilon^2 U(\epsilon(\alpha-t),\epsilon^3 t) 
              + \epsilon^2 V(\epsilon(\alpha+t),\epsilon^3 t) + O(\epsilon^{4}) 
\end{equation}
where $U$ and $V$ satisfy the KdV equations
\begin{equation}
\label{kdv}
\begin{split}
- 2 \partial_T U &= \third \partial_{\beta_-}^3 U + \frac{3}{2} \partial_{\beta_-}(U^2)\\
  2 \partial_T V &= \third \partial_{\beta_+}^3 V + \frac{3}{2} \partial_{\beta_+}(V^2).
\end{split}\tag{KdV}
\end{equation}
See Figure \ref{main result cartoon}.
Here, $\beta_\pm = \epsilon(\beta \pm t)$ represent long wavelength
moving reference frames and 
$T = \epsilon^3 t$ is the very long time scale coordinate.
For technical reasons (which we discuss later), $-x_\alpha$ is the natural
variable to estimate for the water wave equation.  To lowest order,
$-x_\alpha$ is proportional to the height of the wave.  At higher
order, this ceases to be true, though for  purposes of intuition one can think
of $-x_\alpha$ as representing the wave amplitude.

\begin{figure}
\begin{center}
  \psfragscanon
  \includegraphics{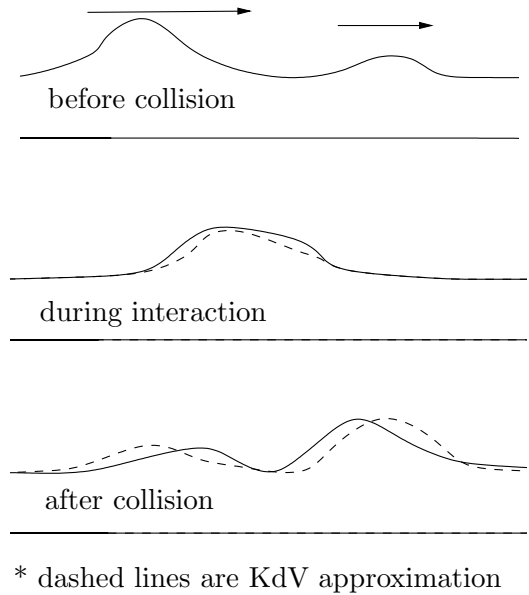}
  \caption{Sketch of the overtaking wave interaction.}
\label{overtaking cartoon}
\end{center}
\end{figure}

The KdV equation was initially derived from the water wave equation
in an attempt to prove the existence of a solitary wave solution
for waves in a canal.  Famously, the KdV equation admits
solitary wave solutions and also multi-soliton solutions.  
(See Figure \ref{overtaking cartoon}).
We will frequently refer to multi-soliton solutions as ``overtaking wave''
collisions.  We remind the reader that the only notable
first order effect after such a collision is that the waves are phase shifted after the collision.   

\begin{figure}
\begin{center}
  \psfragscanon
  \includegraphics{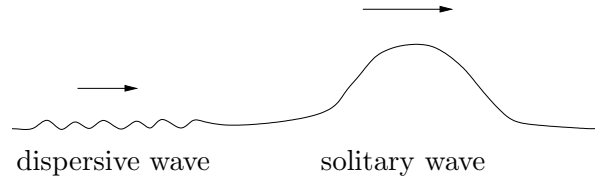}
  \caption{Sketch of the dispersive wave.}
\label{dispersive wave cartoon}
\end{center}
\end{figure}

Given the results in \cite{schneider.etal:00},
one expects to see similar behavior in solutions 
systems modeled by KdV equations.
Though it is unknown if these soliton-like solutions
persist globally, analogous behavior is indeed observed for 
very long times (see \cite{haragus.etal:03}).  
The most notable deviation between true solutions and the KdV approximation
is the size of the phase shift after a collision.  In addition, soliton-like
solutions to the type of systems we study 
frequently develop a very small amplitude
dispersive wave train behind each soliton, which moves in the same direction,
see Figure \ref{dispersive wave cartoon}.    
The KdV approximation does not predict the existence of these dispersive
wave trains.  As these sorts
of discrepancies are observed even in the case where there is only one wave train
moving unidirectionally, we believe that they are, loosely, 
independent of interactions between the left and 
right moving wave trains.  They
reflect intrinsic differences between the approximation and the original system.

\begin{figure}
\begin{center}
  \psfragscanon
  \includegraphics{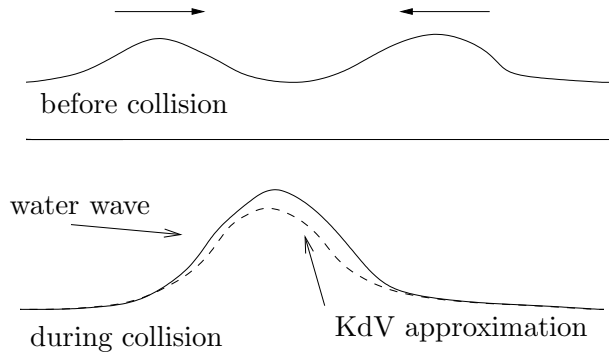}
  \caption{Sketch of the head-on collision.}
\label{head-on cartoon}
\end{center}
\end{figure}

On the other hand, there is evidence that a noticeable interaction
takes place between the left and right moving waves.   One can see from
the form of the approximation in equation (\ref{kdv approx for 
water wave}) that during a head-on 
collision of waves moving
in opposite directions, the KdV approximation predicts that 
the heights of the waves add
linearly.  In true head-on collisions in solutions
to the water wave equation, however,
the height of the waves is slightly different from the sum of the heights
of the waves taken separately---it is slightly larger.  (See the works 
of Maxworthy \cite{maxworthy:76}, Byatt-Smith \cite{byatt-smith:88} \cite{byatt-smith:89},
Cooker, Weidman and Bale \cite{cooker.etal:97}, and Su and Mirie \cite{su.etal:80} \cite{su.etal:82}.)
We sketch this in Figure \ref{head-on cartoon}.

Thus we might expect two types of corrections to the KdV approximation:
\begin{itemize}
\item corrections due to the fact that, even in the case of a purely right (or left)
moving wave train, solutions to the water wave equation are not
exactly described by solutions to the KdV equation.  We will refer to this source of error
as {\it unidirectional error}.  
\item corrections due the fact that the left and right moving wave
trains will interact at higher order.  We call such errors {\it counterpropagation error}.
\end{itemize}
Both of these types of corrections are apparent in our results, and 
to incorporate these two types of corrections, we add an additional three
functions to the KdV wave trains.  
The first two, $F$ and $G$, will correct for unidirectional errors.
The third, $P$, will correct
for counterpropagation errors.
We scale the amplitudes of these three functions so they are
$O(\epsilon^4)$, which is the same as the order of the error in using only the
KdV equations.  $F$ and $G$ will take the same functional 
form as $U$ and $V$, as they correct for differences between 
the approximate and actual wave trains.  So we add  
$$
\epsilon^4 F(\epsilon(\alpha-t),\epsilon^3 t) 
              + \epsilon^4 G(\epsilon(\alpha+t),\epsilon^3 t)
$$
to the first order approximation (\ref{kdv approx for water wave}).  

We do not expect
$P$ to be moving strictly left or right, as it corrects for
the interaction between waves moving in opposite directions.  
So its spatial dependence will be on
$\beta = \epsilon \alpha$.  Suppose that the functions $U$ and $V$ are
solitary wave solutions. In the long wavelength variables we are 
considering, this means
that the right and left moving wave packets are large only over a length
of $O(\epsilon^{-1})$.  In addition, the reference frame moves with
unit velocity.  Thus we expect any interaction of the two waves to last
a time of $O(\epsilon^{-1})$.  Accordingly, we let $P$ depend on
the time variable $\tau = \ep t$.  
That is we add a correction term of the form
$$
\epsilon^4 P(\epsilon \alpha, \epsilon t)
$$
to the KdV model.

Through formal means we find that $P$ satisfies an inhomogeneous wave equation
\begin{equation}
\label{inhomogeneous wave equation}
\begin{split}
\partial_\tau^2 P - \partial_\beta^2 P 
= 3 \partial_\beta^2 \left( U(\beta-\tau,\epsilon^2 \tau) V(\beta+\tau,\epsilon^2 \tau)\right).
\end{split}\tag{IW}
\end{equation}
Similarly, $F$ and $G$ satisfy a pair of driven, linearized KdV equations
\begin{equation}
\label{linearized kdv}
\begin{split}
-2 \pT F = &\third \pbm^3 F + 3 \pbm(UF) + J^-\\
 2 \pT G = &\third \pbp^3 G + 3 \pbp(VG) + J^+.
\end{split}\tag{LK}
\end{equation}
Notice that these equations are linearized about the KdV solutions $U$ and $V$.  The
inhomogeneous terms $J^-$ and $J^+$ are made up of a combination of sums and products 
of $U$, $V$ and $P$.  For explicit forms of these driving terms see the equations
(\ref{lin kdv inhoms water wave}).  Linearized KdV equations are
explicitly solvable, though this is a complicated matter (see
\cite{sachs:83} and \cite{haragus.etal:98}). However, solutions are
simple
to compute numerically.

One can solve inhomogeneous wave equations explicitly 
and easily {\it via} the method of characteristics.  Moreover, we can reduce such
systems to a pair of transport equations by the following fact:  
\begin{fact}
\label{wave to transport}
If $\partial_t f - \partial_x f = 1/2 \partial_x h$ and 
$\partial_t g + \partial_x g = -1/2 \partial_x h$, then $q=f+g$ satisfies
$\partial_t^2 q - \partial_x^2 q = \partial_x^2 h$.  
\end{fact}
Thus, we have
\begin{equation}
\label{transport}
\begin{split}
P^-_\tau + P^-_\beta &= -\frac{3}{2}\partial_\beta(U(\beta-\tau,\ep^2 \tau) V(\beta+\tau,\ep^2 \tau))\\
P^+_\tau - P^+_\beta &=  \frac{3}{2}\partial_\beta(U(\beta-\tau,\ep^2 \tau) V(\beta+\tau,\ep^2 \tau))
\end{split}\tag{T}
\end{equation}
where 
$$
P(\beta,\tau)=P^+(\beta,\tau)+P^-(\beta,\tau).
$$

We remark that the initial data for the modulation equations is determined from initial conditions for the
original system in ways described in Section \ref{error estimates}.  Also, this hierarchy of 
higher order modulation equations is nearly identical to that derived in Wayne and Wright \cite{wayne.etal:02}
for the Boussinesq equation---the chief difference lying in the specific forms of the inhomogeneous terms $J^\pm$.

To enforce the notion of spatial localization,
we will be considering initial data which is of rapid decay, that is,
initial data in
$$
H^s(m) = \left\{ f(\al) | (1 + \al^2)^{m/2} f(\al) \in H^s\right\}.
$$  
The inner product on $H^s(m)$ is given by 
$$
\left( f(\cdot),g(\cdot)\right)_{H^s(m)} = 
\left( (1+\cdot^2)^{m/2} f(\cdot),(1+\cdot^2)^{m/2} g(\cdot) \right)_{H^s},
$$
where we use the standard inner product in $H^s$.  In particular,
the known soliton solutions of the KdV equations
are in such spaces.  

That the KdV equations have solutions for all times with this sort of initial
data is well known.  In particular, we have from \cite{schneider.etal:00}
\begin{theorem}
\label{KdV-solutions} 
Let $\sigma\ge4$.  Then for all $C_I, T_0>0$ there exists $C_1>0$ such that if
$U$, $V$ satisfy (\ref{kdv}) with
initial conditions $U_0$, $V_0$ and 
\begin{equation}
\label{ic kdv bound}
\max \{ \|U_0\|_{H^\sigma(4)\cap H^{\sigma+4}(2) \cap H^{\sigma+9}},
        \|V_0\|_{H^\sigma(4)\cap H^{\sigma+4}(2) \cap H^{\sigma+9}}\} < C_I
\end{equation}
then 
\begin{equation}
\sup_{T\in[0,T_0]} 
\left\{ \|U(\cdot,T)\|_{H^\sigma(4)\cap H^{\sigma+4}(2) \cap H^{\sigma+8}},
        \|V(\cdot,T)\|_{H^\sigma(4)\cap H^{\sigma+4}(2) \cap H^{\sigma+8}}\right\}< C_1.
\end{equation}
\end{theorem}
On the other hand it is less clear that solutions of (\ref{inhomogeneous wave equation})
and (\ref{linearized kdv})
will remain bounded over the very long time scales necessary for the KdV
approximation.  
In \cite{wayne.etal:02} we proved the following result
which guarantees that the solutions of the modulation equations
remain bounded for sufficienly long times.

\begin{proposition}
\label{mods behave}
Fix $T_0 > 0$ and $\sigma>11/2$.  
Suppose, $U_0, V_0$ satisfy (\ref{ic kdv bound}) 
and $U$, $V$, $P^\pm$, $F$ and $G$
satisfy (\ref{kdv}), (\ref{transport}) 
and (\ref{linearized kdv}), 
then there exists 
a constant $C_2$, independent of $\epsilon$, 
such that the solutions of (\ref{transport}) and 
(\ref{linearized kdv}) satisfy the
estimates below:
$$
\sup_{\tau\in[0,T_0\epsilon^{-2}]} \|P^\pm(\cdot,\tau)\|_{H^{\sigma+3}} \le C_2
$$
$$
\sup_{T \in [0,T_0]} \left\{ \|F(\cdot,T)\|_{H^{\sigma} \cap H^{\sigma-4}(2)}, \|G(\cdot,T)\|_{H^{\sigma} \cap H^{\sigma-4}(2)}  
\right\} \le C_2.
$$
Moreover,
$P^\pm(\beta,\tau) = \varphi^\pm(\bpm,T)$
with
$$
\sup_{T \in [0,T_0]} \|\varphi^\pm(\cdot,T)\|_{H^{\sigma+3}\cap H^{\sigma-1}(2)} \le C_2.
$$
\end{proposition}

Finally we note that since 
$\partial_\tau P = \partial_\beta P^+ - \partial_\beta P^-$ we have
$\|\partial_\tau P\|_s \le \|\partial_\beta P\|_s$.

With this preliminary result in hand we can state our principal results.  
Denote the sum of the modulation functions, properly scaled, as
\begin{equation}
\label{psi d 1}
\begin{split}
-\psi^d\left(\al, t\right)
 =&\phantom{+}\epsilon^2 U\left(\epsilon\left(\al-t\right),\epsilon^3 t\right) + \ep^2 V\left(\epsilon\left(\al+t\right),\epsilon^3 t\right) \\
  &+\epsilon^4 F\left(\epsilon\left(\al-t\right),\epsilon^3 t\right) +\ep^4 G\left(\epsilon\left(\al+t\right),\epsilon^3 t\right)\\ 
  &+\epsilon^4 P\left(\epsilon \al, \epsilon t\right). 
\end{split}
\end{equation}
As mentioned earlier, $(x,y)$ are not
the natural coordinates to study solutions to (\ref{water wave}).  The coordinates
we use are $(x_\al, y, x_t)^{tr}$.  $x_\al$ is approximated by $\psi^d$, and the functions
$y$ and $x_t$ are approximated by functions we denote $\psi^y$ and $\psi^u$,
respectively.  They are given by
\begin{equation}
\label{psi y 1}
\begin{split}
\psi^y(\al,t) =&\phantom{+}\ep^2 U(\ep(\al - t),\ep^3 t) + \ep^2 V(\ep(\al + t), \ep^3 t) \\
               &+\ep^4 F(\ep(\al - t),\ep^3 t) + \ep^4 G(\ep(\al + t), \ep^3 t) + \ep^4 P(\ep \al, \ep t) \\
               &+\third \ep^4 \pbm^2 U(\ep(\al - t),\ep^3 t) + \third\ep^4\pbp^2   V(\ep(\al + t), \ep^3 t) \\
               &+\ep^4 \left(U(\ep(\al - t),\ep^3 t) + V(\ep(\al + t), \ep^3 t)\right)^2
\end{split}
\end{equation}
and 
\begin{equation}
\label{psi u 1}
\begin{split}
\psi^u(\al,t) =&\phantom{+}\ep^2 U(\ep(\al - t),\ep^3 t) - \ep^2 V(\ep(\al + t), \ep^3 t) \\
               &+\ep^4 F(\ep(\al - t),\ep^3 t) - \ep^4 G(\ep(\al + t), \ep^3 t) \\
               &+\ep^4 \varphi^{-}(\ep(\al - t),\ep^3 t) - \ep^4 \varphi^{+}(\ep(\al + t), \ep^3 t)\\
               &+\frac{1}{6}\ep^4  \pbm^2 U(\ep(\al - t),\ep^3 t) -  \frac{1}{6}\ep^4 \pbp^2 V(\ep(\al + t), \ep^3 t)\\
               &+\frac{3}{4}\ep^4 U^2(\ep(\al - t),\ep^3 t)- \frac{3}{4}\ep^4  V^2(\ep(\al + t), \ep^3 t).
\end{split}
\end{equation}
We discuss the origin of these equations in Section \ref{derivation}.  

The approximation will be valid in the space
$$
\sob^s = H^s \times H^s \times H^{s-1/2}.
$$
Our main result is:
\begin{theorem}
\label{main result water wave}
Fix $T_0,\ C_I > 0,\ s>4,\ \sigma \ge s+7$.  
Suppose $U$, $V$, $P$, $F$ and $G$ satisfy equations (\ref{kdv}),
(\ref{inhomogeneous wave equation}) and (\ref{linearized kdv}) and
$\psi^d$, $\psi^y$ and $\psi^u$ are the combinations  of these functions
given in (\ref{psi d 1}), (\ref{psi y 1}) and (\ref{psi u 1}).  
Then there exist $\epsilon_0 > 0$ and $C_F > 0$ such that the
following is true.  
If the initial conditions for (\ref{water wave}) are of the form
$\left((x_\alpha(\alpha,0), y(\alpha,0), x_t(\al,0)\right)^{tr} =  
 \left(0, \ep^2 \Theta_y(\ep \alpha), \ep^2 \Theta_u(\ep \al)\right)^{tr}$ 
with 
$$
\max_{i=y,z} \left\{ \|\Theta_i (\cdot)\|_{H^\sigma(4) \cap H^{\sigma+4}(2) \cap H^{\sigma+9}} \right\} \le C_I
$$
then for $\epsilon \in (0,\epsilon_0)$, there is a reparameterization
of the free surface such that the unique solution 
to (\ref{water wave}) satisfies
$$
\left \|
\left( \begin{array}{c}
x_\alpha(\cdot,t)\\
 y(\cdot,t)\\
 x_t(\cdot,t)
\end{array}\right) 
-
\left( \begin{array}{c}
\psi^d (\cdot,t)\\
\psi^y (\cdot,t)\\
\psi^u (\cdot,t)
\end{array}\right) 
\right \|_{\sob^s}
\le C_F \epsilon^{11/2}
$$
for $t \in [0,T_0 \epsilon^{-3}]$.  The constant $C_F$ does not depend on $\epsilon$.
\end{theorem}

\begin{remark}  The loss of the one half power of $\epsilon$ in Theorem \ref{main result water wave} is caused
by the long wave scaling and not a lack of sharpness in the estimates.  
\end{remark}

\begin{remark}  It is clear that the form of the initial conditions specified in the hypotheses
of this theorem do not agree with those found by setting $t = 0$ in the approximation
inequality (unless, of course, $\psi^d \equiv 0$).  This is precisely why we mention the need
to reparameterize the free surface.  We discuss this at length in Section \ref{error estimates}.
\end{remark}

Less technically, this theorem states that  
solutions to (\ref{water wave}),
in the long wavelength limit,
satisfy
$$
x_\al (\al,t)=\psi^d(\al,t) + O(\ep^6)
$$
for times of $O(\epsilon^{-3})$.  
This is, as expected, a marked improvement over the use of KdV alone.

We note that this is not the first time that linearized KdV equations
have been put forward as a means to improve the accuracy of the KdV
approximation. 
Other instances where linearized KdV equations appear include Sachs
\cite{sachs:84}, Sattinger, Haragus and
Nicholls \cite{haragus.etal:03}, Kodama and Taniuti
\cite{kodama.etal:79} and Drazin \cite{drazin:89}.  Moreover, there
have been numerous models put forward over the years which model water
waves in the same scaling regime we are considering.  We refer the
reader to Kodama \cite{kodama:85}, Olver \cite{olver:83}, Bona,
Pritchard and Scott \cite{bona.etal:81}, Craig and Groves
\cite{craig.etal:94}, Dullin, Gottwald and Holm \cite{dullin.etal:03}, and Bona and Chen
\cite{chen.etal:98}.  Much of the work done in the above papers
pertains to analyzing the behavior of the model equations and not to their
connection to the original system.  A notable exception is the recent
work by Bona, Lannes and Colin \cite{bona.etal:03} wherein they prove
the rigorous validity of a large number of Boussinesq style models.
Our particular combination
of linearized KdV equations with an inhomogeneous wave equation
appears to be unique and is asymptotically the most accurate model for long
wavelength solutions to the water wave equation which has currently
been justified rigorously.

The remainder of this paper is organized as follows.  First, in Section \ref{preliminaries},
we conduct a prelimary discussion of the water wave equation.  
Sections \ref{operator K 1} and \ref{operator K 2} contain a thorough discussion
of the operator $K(x,y)$.  
Then, in Section \ref{derivation} we
derive the higher order modulation equations and prove an important estimate.
In Sections \ref{error estimates}
we prove the validity of the approximation {\it i.e.} 
Theorem \ref{main result water wave}.
Finally, Section \ref{proofs} contains the details for a number of proofs.

\vspace{11pt}
\noindent {\bf Acknowledgments:}  
The NSF generously supported this research under grant DMS-0103915.
As this work was done principally as the author's dissertation, 
special thanks also should go to Gene Wayne, who was the author's
thesis advisor, and to the Department of Mathematics and Statistics at Boston
University, where he attended graduate school.

\section{Preliminaries}
\label{preliminaries}

\begin{figure}
\begin{center}
  \psfragscanon
  \includegraphics{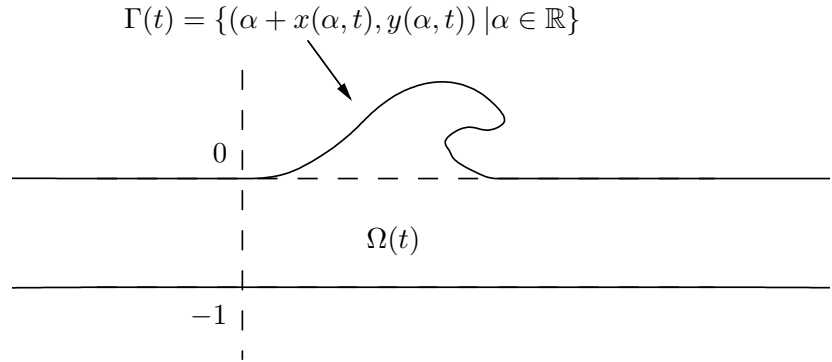}
  \caption{The water wave in Lagrangian coordinates.}
\label{wave set up}
\end{center}
\end{figure}

We begin by discussing the water wave problem in greater detail.
Consider an infinitely long canal of unit mean depth in 
two-dimensions (see Figure \ref{wave set up}).
We denote the region occupied by the fluid at time $t$ as $\Omega(t)$, and
the upper surface as $\Gamma(t)$.  
We parameterize $\Gamma(t)$ by
$\left(\tilde{x}(\al,t), y(\al,t) \right)$, 
where $\al \in \mathbb{R}$ is the parameter, and  $\tilde{x}$ and $y$ are the
real-valued coordinate functions.  It is useful to break $\tilde{x}$ up as follows,
$$
\tilde{x}(\al,t) = \al + x(\al,t).
$$
We consider fluids which are inviscid and incompressible and 
flows which are irrotational.  Also,
we assume that the pressure on the top surface is a constant, and that the
acceleration due to gravity is $1$.  With these assumptions,
the evolution of $x$ and $y$ are given
by the equations
\begin{equation}
\begin{split}
x_{tt}(1+x_\alpha) &+ y_\alpha(1+y_{tt}) = 0\\
y_t &= K(x,y) x_t.
\end{split}\tag{WW}
\end{equation}
(See \cite{craig:85}.)

The first of these two equations is found from Euler's equations for fluid motion.  The
operator $K$ in the second line is a transformation which is linear 
in $x_t$, but depends nonlinearly on $(x,y)$.
That such an operator exists and gives a relationship between 
$x_t$ and $y_t$ is discussed in Sections 
\ref{operator K 1} and \ref{operator K 2}, along with an analysis of $K$.  Much of the difficulty in
answering questions about the water wave equation is related to this operator.  
If the surface of the water is perfectly flat, {\it i.e.} $x=y=0$, then we have
$K(0,0) = K_0$ where $K_0$ is a linear operator defined by 
$\widehat{K_0 f}(k) = \widehat{K_0}(k) \widehat{f}(k)$, with
$\widehat{K_0}(k) = -i \tanh(k)$.  Notice that since $\tanh$ is
a bounded function, $K_0$ is bounded from $H^s$ to $H^s$.  
$K_0$ will be appearing frequently.  

This formulation of the water wave problem is said to be in {\it Lagrangian}, or 
{\it material}, coordinates.  In this point of view we are not in a fixed ``lab''
frame, but instead we are tracking the position of each ``particle'' of water
separately.  That is, $(\al+x(\al,t),y(\al,t))$ gives the location of the
particle which was initially at $(\al+x(\al,0),y(\al,0))$.  The {\it laboratory},
or {\it Eulerian}, point of view is to fix a system of coordinates on the 
fluid domain and to measure the velocity of the fluid at each point of
this fixed reference frame.  For our purposes, it is far more convenient
to use Lagrangian coordinates, however experimentalists work with Eulerian
coordinates.  
We will give formulae for the
approximation in terms of Eulerian coordinates in a future
publication.  The interested reader may also see the author's thesis \cite{wright:03}
for this information.

Since the water wave equation is second order in time for both $x$ and $y$,
one might suppose that four functions are neccessary to specify the initial state
of the system---$x(\al,0)$, $y(\al,0)$, $x_t(\al,0)$ and $y_t(\al,0)$.  In fact,
in general only three are needed, and for initial data with
small amplitudes, only two.  The relationship $y_t = K(x,y) x_t$ 
specifies the value of $y_t(\al,0)$ given the other three functions.  We can
also ``do away'' with the initial condition for $x$, provided we
are in the small amplitude, long wave limit.  If $x(\al,0)$ and its
first derivative are sufficiently small, then $\al + x(\al,0)$ will be invertible.  
This implies that $\Gamma(0)$ will be
a graph over the horizontal coordinate.  And so, without loss of generality,
we can reparameterize the initial conditions so that
$$
\Gamma(0) = \left\{ (\al,y(\al,0)) | \al \in \mathbb{R} \right\}.
$$
Thus we need only to choose $y(\al,0)$ and $x_t(\al,0)$.  As it turns out, we need
to reparameterize the system one more time to prove the approximation theorem, but
we will leave this technicality until Section \ref{error estimates}.  The essential point here
is that due to the freedom in choosing the initial parameterization, we can eliminate
two of our initial conditions.

Even though we can assume that $x(\al,0) = 0$ (or, alternately, is
small), this coordinate grows linearly in time.  
(See the linear
estimates in Chapter 2 of \cite{schneider.etal:00}.)  
As we are concerned
with very long time scales, this is a problem.  As was shown in \cite{schneider.etal:00},
one can replace $x$ with the new coordinate $z=K_0 x$, which is
well-behaved over long times.  Rewriting (\ref{water wave})
with this new variable, we have
\begin{equation}
\label{three-d water wave}
\begin{split}
\partial_t z &= K_0 u\\
\partial_t y &= K(z,y) u\\
\partial_t u &= - \frac{\partial_\al y (1 + \partial_t^2 y)}
                       {1 + K_0^{-1} \partial_\al z}.
\end{split}\tag{WW3}
\end{equation}
We will see that the operator $K(x,y)$ in truth
depends not on $x$ but on $z$, so the abuse of the notation $K(z,y)$ above
is in some sense legitimate (see Section \ref{operator K 1} for further
discussion).  Furthermore, even though $K_0^{-1}$ is not
well-defined, as it blows up at frequency $k=0$, the composition
$$L = K_0^{-1} \partial_\al$$ 
is well-defined as its symbol, $-k/\tanh(k)$,
has no singularities; L is also invertible.  Finally, we notice that the
Maclaurin expansion of $\widehat {K_0}(k) = -ik + O(k^3)$.  Thus
to lowest order
$K_0 \sim -\partial_\alpha$.   And so $z \sim -x_\alpha$.  This is precisely
the reason why, in the Introduction, 
we stated that $-x_\alpha$ is a natural coordinate
for the water wave equation.

Though we will primarily be working with the three-dimensional system
\ref{three-d water wave},
we will need to embed this system into a four-dimensional
system to prove certain aspects of the validity
of the approximation.   We introduce the new coordinate $a = u_t$, and
(\ref{three-d water wave}) becomes
\begin{equation}
\label{four-d water wave}
\begin{split}
\partial_t z &= K_0 u\\
\partial_t y &= K(z,y) u\\
\partial_t u &= a\\
\partial_t a &= 
-\frac{a \partial_\alpha u + \partial_\alpha \partial_t(1 + \partial_t^2 y) 
                           + \partial_t y \partial_t^3 y}
                      {1 + K_0^{-1} \partial_\alpha z}
\end{split}\tag{WW4}
\end{equation}
Though things appear to be getting out of hand, we remark that this is
as large a system as we will need.
Results in \cite{schneider.etal:00}
prove that solutions to equations (\ref{water wave}), (\ref{three-d water wave})
and (\ref{four-d water wave}) do indeed exist for long times.  
We will be considering solutions to the four dimensional system which are 
in 
$$
\sob^s_e = H^s \times H^s \times H^{s-1/2} \times H^{s-1}.
$$

The main goal of this paper is to prove Theorem \ref{main result water wave}.  To do
this, we first prove a similar theorem for solutions to 
(\ref{three-d water wave}), from which Theorem \ref{main result water wave} will follow.
Let 
\begin{equation}
\label{psi d}
\begin{split}
\pd(\al,t) = \psi^d(\al,t) + \ep^6 W_3(\ep \al, \ep t).
\end{split}
\end{equation}
The additional function
$W_3$ will solve an equation we specify later.  While we assure the reader that we will be following Chekhov's rule
and that this gun which appears in the first act will be fired in the third, interested parties
may look ahead to equation (\ref{gun}) in Section \ref{derivation} for more information about $W_3$.
Define the functions
\begin{equation*}
\begin{split}
\pz = & L^{-1} \pd \\
\py = & \pz + \ep^4 \Delta_1 + \ep^6 \Delta_2\\
\pu = & \pa^{-1} \pt \pd.
\end{split}
\end{equation*}
$\Delta_1$ and $\Delta_2$ are combinations of solutions to the modulation equations and are
given in equations ($\ref{delta 1}$) and ($\ref{delta 2}$) in Section \ref{derivation}.  We
justify the presense of the inverse derivative in $\pu$ by means of example.  $\pd$ contains
the term $U$ which solves a KdV equation.  Thus $\pt \pd$ contains the terms from the right hand side
of the KdV equation, all of which are perfect space derivatives.  These then are cancelled by the inverse
derivative.  It is simple to check that this method applies to all terms in $\pd$.

With these function, we have
\begin{theorem}
\label{main result water wave 2}
Fix $T_0,\ C_I > 0,\ s>4,\ \sigma \ge s+7$.  
Let $\pd$, $\pz$, $\py$ and $\pu$ be as above.
Moreover assume the initial conditions for (\ref{kdv})
satisfy
\begin{equation*}
\max \{ \|U_0\|_{H^\sigma(4)\cap H^{\sigma+4}(2) \cap H^{\sigma+9}},
        \|V_0\|_{H^\sigma(4)\cap H^{\sigma+4}(2) \cap H^{\sigma+9}}\} < C_I.
\end{equation*}
Then there exist $\epsilon_0 > 0$ and $C_F=C_F(T_0,C_I,s) > 0$ such that if
the initial conditions for (\ref{three-d water wave}) are of the form
$$
\left( \begin{array}{c}
 z(\al,0)\\
 y(\al,0)\\
 u(\al,0)
\end{array}\right) 
=
\left( \begin{array}{c}
\pz(\al,0)\\ 
\py(\al,0)\\ 
\pu(\al,0) 
\end{array}\right) 
+ \ep^{11/2}\bar{R}_0(\al)
$$
with $\|\bar{R}_0\|_{\sob^s} \le C_I$
then the unique solution of (\ref{three-d water wave}) satisfies the estimate
$$
\left \|
\left( \begin{array}{c}
 z(\cdot,t)\\
 y(\cdot,t)\\
 u(\cdot,t)
\end{array}\right) 
-
\left( \begin{array}{c}
 \pz (\cdot,t)\\
 \py (\cdot,t)\\
 \pu (\cdot,t)
\end{array}\right) 
\right \|_{\sob^s}
\le C_F \epsilon^{11/2}
$$
for $t \in [0,T_0 \epsilon^{-3}]$.  The constant $C_F$ does not depend on $\epsilon$.
\end{theorem}

\section{The operator $K(x,y)$ Part I:  Basics and basic expansions.}
\label{operator K 1}

This is the first of two sections where we discuss the operator $K(x,y)$ which gives the relation between
$y_t$ and $x_t$ in the water wave equation.  Here we briefly discuss the origin of
this operator and 
report expansions of $K$ found in previous work.   We also state
some very basic facts about these expansions.  In Section \ref{operator K 2}
we quote more complicated results and prove some new technical extensions
needed for our purposes.

$K(x,y)$ (or, more precisely $-K(x,y)$)
is sometimes called the Hilbert transform for the region $\Omega(t)$.  
Loosely, given a region $\Omega$ in the complex plane, and any function
$F$ which is analytic in $\Omega$, 
the Hilbert transform for $\Omega$, $H(\Omega)$, is 
a linear operator which relates the real and imaginary parts of 
$F$ on the boundary of $\Omega$.   That is
$$
\textrm{Im} (F) |_{\partial \Omega} = H(\Omega) \textrm{Re} (F) |_{\partial \Omega}.
$$
For example if $\Omega$ were the lower half-plane, 
then the Hilbert transform would be the operator  
$H$, given by $\widehat{H f} = i\ \textrm{sgn}(k) \widehat{f}$. (This particular
operator $H$  is also frequently called {\it the} Hilbert transform.)  
The nature of the operator depends greatly on the region begin studied.  
Unsurprisingly, the proof that such an operator exists is connected to 
the Riemann mapping theorem and to techniques for solving boundary value
problems for Laplace's equation in the plane.  
In this problem, since the region $\Omega(t)$ is completely specified by the coordinate
functions $x$ and $y$, we denote the Hilbert transform by $H(\Omega(t))=-K(x,y)$.  

As we are considering a fluid which is incompressible and a flow which is irrotational, 
$x_t(\al,t) - i y_t(\al,t)$ is the value 
on the upper boundary $\Gamma(t)$, of an analytic function on $\Omega(t)$, 
$\omega(\al + i \beta) = v^x(\al,\beta) - i v^y(\al,\beta)$.  
Here $(v^x,v^y)$ is the velocity field for the fluid in the whole region. 
Thus, given that $K(x,y)$
exists, we have
$$
y_t = K(x,y) x_t.
$$
Of course, the boundary of $\Omega(t)$ is not just $\Gamma(t)$, 
but also includes the bottom of the canal ({\it i.e.} where $\beta = -1$).
As we do not have fluid flow through the bottom, we have $v^y(\al,-1)=0$. 

Under these conditions, $K(x,y)$ has been analyzed extensively by Craig 
in \cite{craig:85} and Schneider and Wayne in \cite{schneider.etal:00}.  
In particular Craig shows that $K(x,y)$ has the following expansion:

\begin{equation}
K(x,y) u  = K_0 u + K_1(x,y) u + S_2(x,y) u 
\end{equation}
where
$$
\widehat{K_0 u} (k) = - i \tanh(k) \widehat{u}(k),
$$
\begin{equation*}
K_1(x,y) u = [x,K_0] \partial_\alpha u - (y + K_0(y K_0)) \partial_\alpha u
\end{equation*}
and $S_2$ is quadratic in $(x,y)$.  

First of all we note that $K_0$ is a bounded operator 
from $H^s$ to $H^s$ since $\tanh(k)$ is a bounded function.  That is
$$
\|K_0 u\|_s \le \|u\|_s.
$$
The operator  
$L=K_0^{-1} \pa$, which is 
well-defined as we discussed in the Section \ref{preliminaries}, will also be used frequently.  
$L$ is not a bounded operator on $H^s$.  It effectively takes one derivative.  That
is
\begin{equation}
\label{L bound}
\|L u\|_s \le \|u\|_{s+1}.
\end{equation} On the other hand, $L^{-1}$ replaces one derivative.  That is,
since 
$$
|\tanh(k)/k| \le (1+k^2)^{-1/2}
$$
we know
\begin{equation}
\label{LI bound}
\|L^{-1} u\|_s \le \|u\|_{s-1}.
\end{equation}

We will be considering functions which are of long wavelength.  That is,
functions of the form $f(\al) = F(\beta)$ where $\beta = \ep \al$.  We define operators
$\Kep$ and $\Lep$ {\it via}
\begin{equation*}
\begin{split}
\Kep F(\beta) = K_0 f(\al)\\
\Lep F(\beta) = L f(\al).
\end{split}
\end{equation*}

Taking the Maclaurin series expansion for $\tanh(k)$ shows that formally
\begin{equation*}
\Kep = - \ep \pb - \frac{1}{3}\ep^3  \pb^3 - \frac{2}{15}\ep^5 \pb^5 + O(\ep^7).  
\end{equation*}
Similarly $\Lep$ and $\Lep^{-1}$ have expansions in terms of derivatives 
\begin{equation*}
\begin{split}
\Lep = &- 1 + \third \ep^2 \pb^2 + \frac{1}{45}\ep^4 \pb^4 + O(\ep^6)\\
\Lep^{-1} =& - 1 - \third \ep^2 \pb^2 - \frac{2}{15}\ep^4 \pb^4 + O(\ep^6).  
\end{split}
\end{equation*}
We call such expansions of Fourier multiplier operators ``long wave approximations''.
The rigorous connection between a long wave
approximation and the original operator is given in the following Lemma, whose
simple proof is contained in Section \ref{proofs}.  
\begin{lemma} 
\label{workhorse}
Suppose $A$ and $A_n$ are linear operators defined by $\widehat{A f}(k)=\widehat{A}(k) \widehat{f}(k)$,
and $\widehat{A_n f}(k)=\widehat{A_n}(k) \widehat{f}(k)$
where $\widehat{A}(k)$ and $\widehat{A_n}(k)$ are complex valued functions.  
Also suppose that 
$|\widehat{A}(k)-\widehat{A_n}(k)| \le C |k|^{n}$. ({\it e.g.} $\widehat{A_n}$ is a Taylor polynomial for
$\widehat{A}$.) Then for $f \in H^{s+n}$  we have
$$
\|A f(\cdot) - A_n f(\cdot)\|_s \le C \|\partial^{n}_\al f(\cdot)\|_{s}.
$$
Moreover, if $f(\al)$ is of long wavelength form---that is if $f(\al)=F(\ep \al)$, 
with $F \in H^{s+n}$---then for $0 < \ep < 1$ there exists
$C$ independent of $\ep$ such that
$$
\|A f(\cdot) - A_n f(\cdot)\|_s \le C \ep^{n-1/2} \|F(\cdot)\|_{s+n}.
$$
\end{lemma}

In \cite{schneider.etal:00}, Schneider and Wayne show that the operator $K(x,y)$ 
does not depend on $x$ {\it per se}, but rather on $z = K_0 x$.  We confuse
the notation for the operators intentionally.  They show the $K(z,y)$ has the following
expansion
$$
K(z,y) u =K_0 u + K_1(z,y) u + S_2(z,y) u. 
$$
where
$$
K_1(z,y) u = M_1(z) \pa u - (y + K_0 y K_0) \pa u
$$
with
\begin{equation*}
\mathfrak{F}[ M_1(z) v ] (k) = 
-\int \frac{ \widehat{K_0}(k) - \widehat{K_0}(l) }{\widehat{K_0}(k-l)} \widehat{z}(k-l) \widehat{v}(l) dl.
\end{equation*}
$S_2$ is an operator which
depends quadratically on $z$ and $y$.  Section \ref{operator K 2} contains an analysis of these
operators. 

By using the hyperbolic trigonometric identity
\begin{equation}
\label{trig id}
\frac{\tanh(l) - \tanh(k)}{\tanh(l-k)} = 1 - \tanh(k) \tanh(l)
\end{equation}
we can simplify the expression for $K_1$ to
\begin{equation}
\label{new K1}
\begin{split}
K_1(z,y) u =& M_1(z+y) \pa u\\
           =& -\left( z+y + K_0(z+y)K_0\right) \pa u
\end{split} 
\end{equation}
Since we know $K_0$ is a bounded operator, it is clear that
\begin{equation}
\label{M1 bound}
\|M_1(z) v \|_s \le \|zv\|_s.
\end{equation}

\section{The derivation}
\label{derivation}

In this section we will derive the higher order correction to the KdV approximation.  
For technical reasons, it is most convenient to work with the water wave equation 
written in the form (\ref{three-d water wave}).
Suppose that one is given the function $\Psi(\al,t)=(\pz(\al,t),\py(\al,t),\pu(\al,t))^{tr}$.  
The amount that this function fails to satisfy (\ref{three-d water wave}) is called 
the {\it residual} and is given by $\resp = (\textrm{Res}_z,\textrm{Res}_y,\textrm{Res}_u)^{tr}$ with
\begin{equation*}
\label{res def 1}
\begin{split}
\textrm{Res}_z&=\partial_t \pz - K_0 \pu\\
\textrm{Res}_y&=\partial_t \py - K(\pz,\py) \pu\\
\textrm{Res}_u&=\pt \pu + 
  \pa \py \frac{1 + \partial_t^2 \py}
              {1 + L \pz}.
\end{split}
\end{equation*}
For a true solution, notice that $\resp$ is identically zero.  

\begin{remark}  We will also consider the four dimensional system (\ref{four-d water wave}).  If
we let $\paa = \pt \pu$, then we have
the additional $\res$ function
\begin{equation*}
\begin{split}
\res_a = & \partial_t \paa +
\frac{\paa \pa  \pu + \pa \partial_t(1 + \partial_t^2 \py) + \partial_t \py \partial_t^3 \py}
                      {1 + L \pz}\\
        = & \partial_t \res_u + \frac{\pa  \pu}
                                     {1 + L \pz} \res_u.
\end{split}
\end{equation*}

\end{remark}

The main goal when deriving modulation equations is to choose a system of equations
such that solutions to this system yield a very small residual.   This is
different than (but connected to) showing that solutions to the modulation equations are close to 
true solutions for the original problem.  This latter issue is precisely that answered
by the main results, Theorem \ref{main result water wave} and Proposition \ref{main result water wave 2},
and is discussed in
Section \ref{error estimates}.  Here, we will perform a series of calculations on the residual
and derive equations (\ref{kdv}), (\ref{inhomogeneous wave equation}) and (\ref{linearized kdv}). In this
process we guarantee the smallness of the residual.  While several of the steps will initially seem to have
little mathematical justification ({\it i.e.} they are formal), once the calculation 
is completed it will be obvious that all steps are valid.  For example, we will take 
\begin{equation}
\label{pu}
\pu = K_0^{-1} \pt \pz.
\end{equation}
With this choice 
$$
\textrm{Res}_z = 0,
$$
which is small indeed!
However, $K_0^{-1}$ is not in general a well-defined operator.  Nonetheless, when we
eventually select $\pz$, $K_0^{-1} \pt \pz$ will make perfect sense.

We are looking for solutions which are small in amplitude and long in wavelength.  So
we let
\begin{equation*}
\begin{split}
\pz(\al,t) = & \ep^2 Z(\beta,\tau) \\
           = &\ep^2 Z_1(\beta,\tau) + \ep^4 Z_2(\beta,\tau) + \ep^6 Z_3(\beta,\tau).
\end{split}
\end{equation*}
Recall $\beta=\ep \al$ and $\tau = \ep t$.  We require $\resp$ to be $O(\ep^{17/2})$. Loosely, 
we need three powers of $\ep$ more than the expected error of $O(\ep^{11/2})$ to account for the long times 
($O(\ep^{-3})$) over which our approximation will be a valid. See Schneider and Wayne \cite{schneider.etal:00} and
Wayne and Wright \cite{wayne.etal:02}.  

\begin{remark} More specifically, if we wish to prove Theorem \ref{main result water wave} in
the space $\mathfrak{H}^s$ we will need
\begin{equation*}
\label{HOPE}
\begin{split}
\|\res_z\|_s \le & C \ep^{17/2}\\
\|\res_y\|_s \le & C \ep^{17/2}\\
\|\res_u\|_{s-1} \le & C \ep^{17/2}\\
\|\res_a\|_{s-1} \le & C \ep^{19/2}
\end{split}
\end{equation*}
for $0 \le t \le T_0 \ep^{-3}$.
Given the definition
of $\res_a$, the final estimate will follow
automatically from the estimate on $\res_u$.  
\end{remark}

We have already chosen $\pu$ in terms of $\pz$.  We will use first use the expression for
$\res_y$ to similarly determine $\py$ in terms of $\pz$.  This is not as simple a matter because while
$K_0$ commutes with $\pt$, the full operator $K(z,y)$ does not.

We have
\begin{equation*}
\begin{split}
\res_y = &-\pt \py + K_0 \pu + M_1(\pz + \py) \pa \pu + S_2(\py,\pz) \pu\\
       = &-\pt \py + \pt \pz + M_1(\pz + \py) K_0^{-1} \pa \pt \pz + S_2(\py,\pz) \pu.
\end{split}
\end{equation*}
Notice that in the above
expression we can cancel the linear terms by taking $\py \sim \pz$.  More precisely, we set
$$
\py(\al,t) = \ep^2 Z(\beta,\tau) + \ep^4 \Delta_1(\beta,\tau) + \ep^6 \Delta_2 (\beta,\tau),
$$
for as yet undetermined functions $\Delta_i$.  So
\begin{equation*}
\begin{split}
\res_y = &-\ep^5 \ptau \Delta_1 - \ep^7 \ptau \Delta_2  + M_1(2 \ep^2 Z) \Lep \ep^3 \ptau Z\\
         &+M_1(\ep^4 \Delta_1 + \ep^6 \Delta_2) \Lep \ep^3 \ptau Z + S_2(\py,\pz) \pu\\
       = &-\ep^5 \ptau \Delta_1 - \ep^7 \ptau \Delta_2  - 2 \ep^5 Z \Lep \ptau Z -  \ep^6 \Kep \left( 2 Z \pb \ptau  Z\right) \\
         &-\ep^7 \Delta_1 \Lep \ptau Z - \ep^8 \Kep \left(\Delta_1 \pb \ptau Z\right)  + \ep^9 M_1(\Delta_2) \Lep \ptau Z\\
         & + S_2(\py,\pz) \pu.
\end{split}
\end{equation*}

A number of the terms in $\res_y$
are already $O(\ep^{17/2})$.  
By Lemma \ref{workhorse} we have the following estimate on $\Kep$:
\begin{equation}
\label{Kep small}
\begin{split}
\|\Kep F\|_s \le \ep^{1/2} \|F\|_{s+1}.
\end{split}
\end{equation}
Thus terms containing $\Kep$ can be considered to be a power of $\ep$ smaller
than they appear (though this costs a derivative). 
On the other hand,
$K_0$ is bounded so we have
\begin{equation}
\label{Kep bound}
\|\Kep F\|_s \le \ep^{-1/2} \|F\|_s.
\end{equation}
Thus we can use $\Kep$ either as a bounded functional {\it or} to gain powers of $\ep$,
but not both.  Notice that $\Lep$ does not contribute any additional powers of $\ep$ in any case.  
We separate out all the terms that are already sufficiently small into error terms.  That is
\begin{equation*}
\begin{split}
\res_y = &-\ep^5 \ptau \Delta_1 - \ep^7 \ptau \Delta_2  - 2 \ep^5 Z_1 \Lep \ptau Z_1 \\
         &- 2\ep^7 Z_1 \Lep \ptau Z_2 - 2 \ep^7 Z_2 \Lep \ptau Z_1  \\
         &- \ep^6 \Kep\left( 2 Z_1 \pb \ptau Z_1\right) - \ep^7 \Delta_1 \Lep \ptau Z_1\\ &+ E^y_{small} + E^y_{S_2}
\end{split}
\end{equation*}
with
\begin{equation*}
\begin{split}
E^y_{small} = & 2 \ep^9 Z_3 \Lep \ptau  Z + 2 \ep^9 (Z_1 + \ep^2 Z_2) \Lep \ptau  Z_3 + 2 \ep^9 Z_2 \Lep \ptau Z_2\\
      &- \ep^8 \Kep \left( (2 Z_2+2\ep^2 Z_3) \pb \ptau  Z\right) -\ep^8 \Kep\left(2 Z_1 \pb\ptau (Z_2 +\ep^2 Z_3)\right)\\
      &-\ep^9 \Delta_1 \Lep \ptau (Z_2 + \ep^2 Z_3) - \ep^8 \Kep \left(\Delta_1 \pb \ptau Z\right)  + \ep^9 M_1(\Delta_2) \Lep \ptau Z\\
E^y_{S_2} = &S_2(\py,\pz) \pu.
\end{split}
\end{equation*}

It is clear that $E^y_{small}$ is $O(\ep^{17/2})$.  That is
$$
\|E^y_{small}\|_s \le C \ep^{17/2}.
$$
The constant $C$ depends on various norms of the functions $Z_i$, $\ptau Z_i$ and $\Delta_i$.  
Specifically, chasing through the various terms in $E^y_{small}$ and applying the estimates
in (\ref{L bound}), (\ref{M1 bound}), (\ref{Kep small}) and (\ref{Kep bound}),
one can show that $C$ depends on
$\|Z_1\|_{s+1}$, $\|Z_2\|_{s+1}$, $\|Z_3\|_{s}$,
$\|\ptau Z_1\|_{s+2}$, $\|\ptau Z_2\|_{s+2}$, $\|\ptau Z_3\|_{s+1}$,
$\|\Delta_1\|_{s+1}$ and $\|\Delta_2\|_{s}$.

The term $E^y_{S_2}$ is also $O(\ep^{17/2})$ though this is not as obvious.  We prove this in Proposition
\ref{S_2 special} in Section \ref{operator K 2}.  The proof of this relies strongly on the fact that we have taken
$\py$ and $\pz$ such that $\py - \pz$ is $O(\ep^4)$.  This causes a cancellation in $S_2$, which in turn
makes this term small.

We now expand $\Lep$ and $\Kep$ in the remaining low order terms in $\res_y$ to find
\begin{equation*}
\begin{split}
\res_y = &-\ep^5 \ptau \Delta_1 - \ep^7 \ptau \Delta_2  + 2 \ep^5 Z_1  \ptau Z_1 + \frac{4}{3} \ep^7 Z_1 \ptau \pb^2  Z_1\\
         &+ 2\ep^7 Z_1  \ptau Z_2 + 2 \ep^7 Z_2  \ptau Z_1  \\
         &+ 2 \ep^7 \pb Z_1 \pb \ptau  Z_1  + \ep^7 \Delta_1  \ptau  Z_1\\& + E^y_{small} + E^y_{S_2}+E^y_{lwa}
\end{split}
\end{equation*}
with
\begin{equation*}
\begin{split}
E^y_{lwa} = &-2 \ep^5 Z_1 (\Lep+1 - \third \ep^2 \pb^2) \ptau  Z_1 - 2\ep^7 Z_1 (\Lep+1) \ptau Z_2\\
       & -2 \ep^7 Z_2 (\Lep+1)\ptau  Z_1
      -2 \ep^6 (\Kep + \ep\pb)( Z_1 \pb \ptau Z_1)\\
       & - \ep^7 \Delta_1 (\Lep+1) \ptau  Z_1.
\end{split}
\end{equation*}
Each term in $E^y_{lwa}$ is $O(\ep^{17/2})$ by Lemma \ref{workhorse}.  That is
$$
\|E^y_{lwa}\| \le C\ep^{17/2}
$$
where $C$ depends on
$\|Z_1\|_{s+3}$,$\|Z_2\|_{s}$,
$\|\ptau Z_1\|_{s+4}$,$\|\ptau Z_2\|_{s+2}$
and $\|\Delta_1\|_{s}$.
(The subscript ``lwa'' stands for ``long wave approximation''.)

The only $O(\ep^5)$ terms remaining in $\res_y$ are
$$
- \ep^5 \ptau \Delta_1 + 2\ep^5 Z_1 \ptau Z_1
$$
which we remove by selecting
\begin{equation}
\label{delta 1}
 \Delta_1 = Z_1^2.
\end{equation}
So
\begin{equation*}
\begin{split}
\res_y = & - \ep^7 \ptau \Delta_2  + \frac{4}{3} \ep^7 Z_1 \ptau \pb^2 Z_1\\
         &+ 2\ep^7 Z_1  \ptau Z_2 + 2 \ep^7 Z_2  \ptau  Z_1  \\
         &+ 2 \ep^7 \pb Z_1 \pb \ptau  Z_1  + \ep^7 Z_1^2  \ptau  Z_1\\
         & + E^y_{small} + E^y_{S_2}+E^y_{lwa}.
\end{split}
\end{equation*}

The remaining $O(\ep^7)$ terms in $\res_y$ are all perfect time derivatives with the exception of
$$
\frac{4}{3}  Z_1 \ptau \pb^2  Z_1.
$$
Notice, however, that
\begin{equation*}
\begin{split}
 &  \ptau \left( \frac{4}{3} Z_1 \pb^2  Z_1  - \frac{2}{3} (\ptau  Z_1)^2 \right) \\
=& \frac{4}{3}  Z_1 \ptau \pb^2 Z_1 + \frac{4}{3} \ptau  Z_2\left(\pb^2 Z_1 - \ptau^2 Z_1\right).
\end{split}
\end{equation*}
Given the form of the approximation in equation (\ref{kdv approx for water wave}), 
it is not unreasonable to suspect that
\begin{equation}
\label{switch}
\pb^2  Z_1 - \ptau^2 Z_1 \sim O(\ep^2).
\end{equation}

We are now in a position to select $\Delta_2$.  
Taking
\begin{equation}
\label{delta 2}
\begin{split}
\Delta_2 = &(\pb Z_1)^2 + 2 Z_1 Z_2 + \third Z_1^3\\ 
    &+ \frac{4}{3}  Z_1 \pb^2  Z_1 - \frac{2}{3}(\ptau Z_1)^2
\end{split}
\end{equation}
gives
$$
\res_y = E^y_{small} + E^y_{S_2} + E^y_{lwa} + E^y_{switch}
$$
where
$$
E^y_{switch} = \frac{4}{3} \ep^7 \ptau Z_1 (\pb^2 Z_1 - \ptau^2  Z_1).
$$
Given that our assumption (\ref{switch}) is valid, we have shown that with our choices
for $\py$ and $\pu$ in terms of $\pz$ that $\res_y = O(\ep^{17/2})$.  More specifically,
we have shown that if
$$
\|\pb^2 Z_1 - \ptau^2  Z_1\|_s \le C\ep^{3/2},
$$
then
$$
\res_y \le C \ep^{17/2}
$$
where the constant depends only on
$\|Z_1\|_{s+3}$, $\|Z_2\|_{s+1}$, $\|Z_3\|_{s}$,
$\|\ptau Z_1\|_{s+4}$, $\|\ptau Z_2\|_{s+2}$, $\|\ptau Z_3\|_{s+1}$.

Now that we have computed $\py$ and $\pu$ in terms of $\pz$, we now turn our attention to determining
$\pz$ by examining $\res_u$.  
\begin{equation*}
\begin{split}
\res_u = K_0^{-1} \pt^2 \pz + \pa\py \frac{1+\pt^2 \py}{1 + L \pz}.
\end{split}
\end{equation*}
We expand $(1+L \pz)^{-1}$ by the geometric series to find
\begin{equation*}
\begin{split}
\res_u = K_0^{-1} \pt^2 \pz + \pa\py (1+\pt^2 \py)(1 - L \pz + (L \pz)^2) + E^u_1
\end{split}
\end{equation*}
where
$$
E^u_{geo} = \pa\py \frac{1+\pt^2 \py}{1 + L \pz} - \pa\py (1+\pt^2 \py)(1 - L \pz + (L \pz)^2).
$$
Since $\pz$ is ``small'' this error term can be shown to be $O(\ep^{17/2})$.
We have
\begin{lemma}
\label{geometric}
Let $f \in H^{s+1}$, $s>0$.  Take $\ep_0$ such that $\ep_0^2 \|f\|_{L^\infty} \le 1/2$.
Then for $0 < \ep < \ep_0$ we have that the function
$$
g(\ep x) = \frac{1}{1+\ep^2 f(\ep x)} - 1 + \ep^2 f(\ep x)
$$
satisfies $\|g(\ep \cdot)\|_s \le C \ep^{7/2}$ for $C$ independent of $\ep$.
\end{lemma}
\begin{remark}
Under the same hypotheses as in Lemma \ref{geometric}, arguments similar to the proof of that Lemma
show:
\begin{itemize}
\item
$
(1+\ep^2 f(\ep x))^{-1} \in C^{s}
$
and is bounded there independent of $\ep$,
\item
$
(1+\ep^2 f(\ep x))^{-1} - 1  \in H^s
$
and has norm there bounded by $C \ep^{3/2}$ for $C$ independent of $\ep$ and 
\item
$
(1+\ep^2 f(\ep x))^{-1} - 1  + \ep^2 f(\ep x) -\ep^4 f^2(\ep x) \in H^s
$
and has norm there bounded by $C \ep^{11/2}$ for $C$ independent of $\ep$.
\end{itemize}
\end{remark}

Now, after substituting in from the definitions of $\py$ and $\pz$, we collect all the terms which are smaller than
$O(\ep^{17/2})$ and find
\begin{equation*}
\begin{split}
\res_u = &\ep^4 \Kep^{-1} \ptau^2 (Z_1 + \ep^2 Z_2 + \ep^4 Z_3) + \ep^3 \pb Z_1 - \ep^5 \pb Z_1 \Lep Z_1\\
         &+\ep^5 \pb\left(Z_2 + Z_1^2\right) + \ep^7 \pb Z_1 \left(\pb^2  Z_1 - \Lep Z_2 + (\Lep Z_1)^2\right)\\
         &-\ep^7 \pb(Z_2 + Z_1^2) \Lep Z_1 + \ep^7 \pb  Z_3\\ 
         &+\ep^7  \pb \left( (\pb  Z_1)^2 + 2 Z_1  Z_2 + \third Z_1^3 \right)\\ 
         &+ \ep^7 \pb \left(\frac{4}{3} Z_1 \pb^2 Z_1 - \frac{2}{3}(\ptau  Z_1)^2\right)\\
         &+ E^u_{geo} + E^u_{small} + E^u_{switch}.
\end{split}
\end{equation*}
We omit the exact expression for $E^u_{small}$ because it is both lengthy and uninteresting. We have
\begin{equation}
\label{troublemaker}
\|E^u_{small}\|_s \le C\ep^{17/2}
\end{equation}
where the constant $C$
depends on $H^{s+1}$ norms of the functions $Z_i$ and the $H^s$ norms of
$\ptau^2 Z_i$.  We have also replaced one instance of $\ptau^2 Z_1$ with $\pb^2  Z_1$ (much
as we did earlier), thus
the term
$$
E^u_{switch} = \ep^7 \pb Z_1 \left(\pb^2 Z_1 - \ptau^2  Z_1\right).
$$

We now define new functions $W_i$ by
$Z_i = \Lep^{-1} W_i$.  This seemingly mysterious (and sudden!) 
change of variables will seem less so if we remind the reader that
at the end of the day we wish model not $z$ but rather the function $x_\al$.  
Accordingly, if we approximate $x_\al$ by a function 
$$
\pd = \ep^2 W_1 + \ep^4 W_2 + \ep^6 W_3,
$$ 
then it is logical
to take
$$
\pz = L^{-1} \pd
$$
and in the long wavelength limit we arrive at these functions $W_i$.

So we have
\begin{equation*}
\begin{split}
\res_u = &\ep^3 \pb^{-1} \ptau^2 ( W_1 + \ep^2 W_2 + \ep^4 W_3) + \ep^2 \Kep W_1 - \ep^4 W_1 \Kep W_1\\
         &+\ep^4 \Kep W_2 + \ep^5 \pb (\Lep^{-1} W_1)^2 \\
&+ \ep^2 \Kep W_1 \left(\ep^3 \Kep \pb W_1 - \ep^4 W_2 + \ep^4  W_1^2\right)\\
         &-(\ep^6 \Kep W_2 +\ep^7 \pb(\Lep^{-1} W_1)^2)W_1 + \ep^6 \Kep W_3 \\
         &+ \ep^5  \pb \left( (\Kep  W_1)^2 + 2 \ep^2 \Lep^{-1} W_1 \Lep^{-1} W_2 + \third \ep^2  (\Lep^{-1} W_1)^3 \right)\\ 
         &+ \ep^7 \pb \left(\frac{4}{3}  \Lep^{-1} W_1 \Lep^{-1} \pb^2 W_1 - \frac{2}{3}(\Lep^{-1} \ptau W_1)^2\right)\\
         &+ E^u_{geo} + E^u_{small} + E^u_{switch}.
\end{split}
\end{equation*}
At this time, the presence of inverse $\beta$ derivatives may seem problematic.  Notice that
each such precedes a time derivative.  Once we select the functions
$W_i$ we will see that there can be an exchange between time and space derivatives, which
will justify the instances of $\pb^{-1}$.

Now we replace $\Kep$  and $\Lep^{-1}$ by their long wave approximates and find
\begin{equation*}
\begin{split}
\res_u = &\ep^3 \pb^{-1} \ptau^2 ( W_1 + \ep^2 W_2 + \ep^4 W_3)\\ 
         &- \ep^2 (\Kt)  W_1\\ 
         &+ \ep^4 W_1 (\Kd)  W_1\\
         &- \ep^4 (\Kd)  W_2 + \ep^5 \pb (W_1)^2 + \frac{2}{3} \ep^7  \pb( W_1 \pb^2 W_1)\\
         & - \ep^7 \pb W_1 \left(- \pb^2 W_1 -  W_2 +   W_1^2\right)\\
         &-\ep^7 W_1 (-\pb W_2 +\pb( W_1)^2) - \ep^7 \pb  W_3 \\
         &+ \ep^7 \pb \left((\pb W_1)^2 + 2  W_1  W_2 - \third  W_1^3 \right)\\ 
         &+ \ep^7 \pb \left(\frac{4}{3}  W_1 \pb^2  W_1 - \frac{2}{3}( \ptau  W_1)^2\right)\\
         &+ E^u_{geo} + E^u_{small} + E^u_{switch} + E^u_{lwa}.
\end{split}
\end{equation*}

The error made by the long wave approximations is denoted $E^u_{lwa}$.  By Lemma \ref{workhorse}
we have
$$
\|E^u_{lwa}\|_s \le C \ep^{17/2}
$$
where $C$ depends on 
$\|W_1\|_{s+7}$, $\|W_2\|_{s+5}$, $\|W_3\|_{s+3}$, $\|\ptau W_1\|_{s+2}$.

Now we organize the above as
\begin{equation}
\begin{split}
\label{master mold}
\res_y = &  \ep^3 \pb^{-1} \ptau^2  W_1 - \ep^3 \pb  W_1  \\
         &-\ep^5 \third \pb^3 W_1 + \ep^5 \frac{3}{2}\pb ( W_1)^2\\
         & +\ep^5 \pb^{-1} \ptau^2  W_2 - \ep^5 \pb  W_2  \\
         & -\ep^7 \third \pb^3 W_2 + \ep^7 3 \pb (W_1  W_2)\\
         &-\frac{2}{15}\ep^7 \pb^5  W_1 + \ep^7 \third  W_1 \pb^3 W_1\\
         &+\ep^7 2\pb( W_1 \pb^2 W_1) + \ep^7 \frac{3}{2}\pb(\pb  W_1)^2\\
         &-\ep^7 \frac{4}{3} \pb ( W_1)^3 - \ep^7 \frac{2}{3} \pb(\ptau  W_1)^2\\
         & +\ep^7 \pb^{-1} \ptau^2 W_3 - \ep^7 \pb  W_3  \\
         &+ E^u_{geo} + E^u_{small} + E^u_{switch} + E^u_{lwa}.
\end{split}
\end{equation}

The term on the first line of 
right hand side looks formally like an inverse derivative of a wave equation:
$$
\pb^{-1} \ptau^2  W_1 - \pb  W_1  = \pb^{-1} \left(\ptau^2 W_1 - \pb^2 W_1\right).
$$
We cancel this term (to lowest order) by taking $W_1$ of the form
$$
W_1(\beta,\tau) = -U(\beta-\tau,\ep^2 \tau) -  V(\beta+\tau,\ep^2 \tau).
$$
Recall $\bpm = \beta \pm \tau$ and $T=\ep^2\tau$.  The ``minus'' signs may seem arbitrary, but
are included at this stage so that they agree with previous work in the area.  Noting that the third line
looks very much like the first, we also set
$$
W_2(\beta,\tau) = -F(\beta-\tau,\ep^2 \tau) -  G(\beta+\tau,\ep^2 \tau) - P(\beta,\tau).
$$
These choices for $W_1$ and $W_2$ are precisely those described heuristically in the Introduction.

The first three lines in (\ref{master mold}) become
\begin{equation*}
\begin{split}
 &\ep^5 \left(2 \pT U + \third \pbm^3 U  + \frac{3}{2}\pbm U^2 \right)\\
+&\ep^5 \left(- 2 \pT V + \third \pbp^3 V + \frac{3}{2}\pbp V^2\right)\\
+&\ep^5 \left( \pb^{-1}\left(\pb^2 P - \ptau^2 P\right) + 3 \pb(UV) \right)\\
+&\ep^7 \left(2 \pT F - 2 \pT G
 - \pb^{-1} \left(\pT^2 U + \pT^2 V\right)\right)\\
-&\ep^9 \pb^{-1} \left(\pT^2 F + \pT^2 G\right).
\end{split}
\end{equation*}
We cancel everything multiplied by $\ep^5$ by taking
\begin{equation*}
\begin{split}
- &2 \pT U = \third \pbm^3 U  + \frac{3}{2}\pbm U^2\\ 
  &2 \pT V = \third \pbp^3 V  + \frac{3}{2}\pbp V^2\\ 
  & \ptau^2 P - \pb^2 P = 3 \pb^2 (UV)
\end{split}
\end{equation*}
which are precisely equations (\ref{kdv}) and (\ref{inhomogeneous wave equation}).  By
Proposition \ref{mods behave} we know the solutions to these equations are well-behaved
over the long time scales.

Given that the functions $U$ and $V$ have been chosen to solve (\ref{kdv}), one
computes that 
\begin{equation*}
\begin{split}
\partial_T^2 U &= \partial_{\beta_-} \left(\frac{1}{36} \partial_{\beta_-}^5 U
+\frac{9}{4} U^2 \partial_{\beta_-} U
+\frac{1}{2} U \partial_{\beta_-}^3 U
+\frac{3}{4} \partial_{\beta_-} U\partial_{\beta_-}^2 U \right)\\
\partial_T^2 V
&= \partial_{\beta_+} \left(\frac{1}{36} \partial_{\beta_+}^5 V
+\frac{9}{4} V^2 \partial_{\beta_+} V
+\frac{1}{2} V\partial_{\beta_+}^3 V
+\frac{3}{4} \partial_{\beta_+} V\partial_{\beta_+}^2 V \right).
\end{split}
\end{equation*}
Thus the term $\pb^{-1} \left(\pT^2 U + \pT^2 V\right)$ is perfectly
well-defined.  For brevity, we will continue to write these terms with the inverse
derivatives instead of in the longer form above.  

Moreover, now we can put more precise estimates on $E^y_{switch}$ and $E^u_{switch}$.  In particular,
since each time derivative for solutions to KdV equations count for three space derivatives,
we have
$$
\|\pb^2 W_1 - \ptau^2 W_1\|_s \le C \ep^{3/2}
$$
where $C$ depends on $\|W_1\|_{s+6}$.

Recall from Fact \ref{wave to transport} and Proposition \ref{mods behave} 
that solutions to (\ref{inhomogeneous wave equation}) can be rewritten as
\begin{equation*}
\begin{split}
P(\beta,\tau) = & P^+(\beta,\tau) + P^-(\beta,\tau)\\
              = & \varphi^+(\bp,T) + \varphi^-(\bm,T).
\end{split}
\end{equation*}
The functions $\varphi^\pm$ are rapidly decaying.  We make this
decomposition so that every remaining term in (\ref{master mold}): 
\begin{itemize}
\item will be a unidirectional term which is rapidly decaying;
\item will be a product of two such terms which are moving in opposite directions;
\item or will include a derivative of $W_3$.  
\end{itemize}
That is
\begin{equation*}
\begin{split}
\res_y = & \ep^7 \left(2 \pT F + \third \pbm^3 F + 3 \pbm(UF) + J^- \right)\\
         & +\ep^7 \left(-2 \pT G + \third \pbp^3 G + 3 \pbp(VG) + J^+ \right)\\
         & +\ep^7 \left(\pb^{-1} \left( \ptau^2 W_3 - \pb^2 W_3\right) + J^s\right) \\
         &-\ep^9 \pb^{-1} \left(\pT^2 F + \pT^2 G\right)\\
         &+ E^u_{geo} + E^u_{small} + E^u_{switch} + E^u_{lwa}+E^u_{time}.
\end{split}
\end{equation*}
where
\begin{equation}
\begin{split}
\label{lin kdv inhoms water wave}
J^- =&\phantom{+}  3 \partial_{\beta_-}(U \varphi^- )
                       +4 U^2 \partial_{\beta_-} U  + \frac{7}{3} U \partial_{\beta_-}^3 U  \\
          &             +\frac{11}{3}  \partial_{\beta_-} U \partial_{\beta_-}^2 U + \frac{2}{15} \pbm^5 U+\third \pbm^3 \varphi^-\\
          &             -\pbm^{-1} \pT^2 U\\ 
J^+ =&\phantom{+}  3 \partial_{\beta_+}(V \varphi^+ )
                       +4 V^2 \partial_{\beta_+} V  + 7/3 V \partial_{\beta_+}^3 V \\ 
          &             +\frac{11}{3} \partial_{\beta_+} V \partial_{\beta_+}^2 V + \frac{2}{15} \pbp^5 V+\third \pbp^3 \varphi^+\\
          &             -\pbp^{-1} \pT^2 V\\ 
J^s =& \phantom{+}\partial_{\beta} \left( U (3 G + 3 \varphi^+ + 4V^2 + \frac{7}{3} \partial_{\beta_+}^2 V)\right)\\
          & +\partial_\beta\left(V (3 F + 3 \varphi^- + 4U^2 + \frac{7}{3} \partial_{\beta_-}^2 U)\right)\\
                                   &+4 \pb \left(\partial_{\beta_-}U\partial_{\beta_+}V\right).
\end{split}
\end{equation}
and
\begin{equation*}
\begin{split}
E^u_{time} = & +\ep^9 \frac{4}{3} \pb \left( \left(\pbm U - \pbp V\right)\left(\pT U + \pT V\right)\right)\\
             & -\ep^{11} \frac{2}{3} \pb \left(\left(\pT U + \pT V\right)^2\right).
\end{split}
\end{equation*}
Notice that $J^\pm=J^\pm(\bpm,T)$.
$E^u_{time}$ (so called because each term in it contains some sort of time derivative) is clearly
$O(\ep^{17/2})$.  That is
$$
\| E^u_{time} \|_s \le C \ep^{17/2}.
$$
The constant above depends on $\|U\|_{s+4}$ and $\|V\|_{s+4}$.  

The term $\ep^9 \pb^{-1} \left(\pT^2 F + \pT^2 G\right)$ is not
included in $E^u_{time}$ for the following reason.  In a moment, when we
select the equations $F$ and $G$ solve, a consequence will be that
there will be terms in $\pT^2 F$ and $\pT^2 G$ which are $O(\ep^{-2})$.

By taking
\begin{equation*}
\begin{split}
-2 \partial_T F&=\third \pbm (U F) + \frac{3}{2}\pbm^3 F + J^-\\
 2 \partial_T G&= \third \pbp (V G) + \frac{3}{2}\pbp^3 G + J^+
 \end{split}
\end{equation*}
we cancel nearly all the terms which are not in the various $E^u$ terms.  These are the
linearized KdV equations (\ref{linearized kdv}) discussed in the Introduction.  Proposition
\ref{mods behave} guarantees that the solutions are well-behaved.
We are left with
\begin{equation*}
\begin{split}
\res_y = & +\ep^7 \left(\pb^{-1}\left( \ptau^2 W_3 - \pb^2 W_3\right) + J^s\right) \\
         &-\ep^9 \pb^{-1} \left(\pT^2 F + \pT^2 G\right)\\
         &+ E^u_{geo} + E^u_{small} + E^u_{switch} + E^u_{lwa}+E^u_{time}.
\end{split}
\end{equation*}

Now we consider the terms in $\pb^{-1} \left(\pT^2 F + \pT^2 G\right)$.  Notice that
\begin{equation*}
\begin{split}
-\partial_{\beta_-}^{-1} \partial_T^2 F = & \partial_{\beta_-}^{-1} \partial_T (\frac{1}{6} \partial_\beta^3 F +\frac{3}{2}\partial_\beta(U F) +\frac{1}{2} J^-)\\
                          =& \frac{1}{6} \pbm^2 \pT F + \frac{3}{2}\pT (U F) + \frac{1}{2} \pbm^{-1} \pT J^-.
\end{split}\end{equation*}
$J^-$ contains the term $3 \pbm (U \varphi^-) + \third \pbm^3 \varphi^-$.  From
the definition of $\varphi^-$ we know
$$
\partial_T \varphi^- = -\ep^{-2} \frac{3}{2} \partial_{\beta} (U V).
$$
So we have
\begin{equation*}
\begin{split}
 &\frac{1}{2} \pbm^{-1} \pT (3 \pbm (U \varphi^-) + \third \pbm^3 \varphi^-)\\
=&\frac{3}{2}\varphi^- \pT U + \frac{3}{2}U \pT \varphi^- 
                         +\frac{1}{6} \pbm^2 \pT \varphi^-\\
=&\frac{3}{2}\varphi^- \pT U - \frac{9}{4}\ep^{-2} U \pb (U V) + \frac{1}{4} \ep^{-2} \pb^3(UV).
\end{split}
\end{equation*}
We treat $\partial_T^2 G$ in the same fashion.  So we can write
\begin{equation*}
\begin{split}
-\ep^9 \pb^{-1} \left(\pT^2 F + \pT^2 G\right) = E^u_{F,G} - \ep^7\left( \frac{9}{4} (U+V) \pb(UV) + \frac{1}{2} \pb^3(UV) \right).
\end{split}
\end{equation*}
By construction $E^u_{F,G}$ satisfies the estimate
$$
\|E^u_{F,G}\|_s \le C \ep^{17/2} 
$$
where $C$ depends on 
$\|U\|_{s+7}$,
$\|V\|_{s+7}$,
$\|F\|_{s+5}$,
and $\|G\|_{s+5}$.

We have
\begin{equation*}
\begin{split}
\res_y = & +\ep^7 \left(\pb^{-1} \left(\ptau^2 W_3 - \pb^2 W_3\right) + J^s\right) \\
         &- \ep^7\left( \frac{9}{4} (U+V) \pb(UV) + \frac{1}{2} \pb^3(UV) \right)\\
         &+ E^u_{geo} + E^u_{small} + E^u_{switch} + E^u_{lwa}+E^u_{time} + E^u_{F,G}.
\end{split}
\end{equation*}
By selecting
\begin{equation}
\label{gun}
\ptau^2 W_3 - \pb^2 W_3 = \pb\left(\frac{9}{4} (U+V) \pb(UV) + \frac{1}{2} \pb^3(UV) - J^s \right) 
\end{equation}
the gun goes off and we cancel all remaining $O(\ep^7)$ terms.  Thus
$$
\res_y =  E^u_{geo} + E^u_{small} + E^u_{switch} + E^u_{lwa}+E^u_{time} + E^u_{F,G}.
$$
Each of the $E^u$ is $O(\ep^{17/2})$.  

Unlike the previous equations (\ref{kdv}), (\ref{inhomogeneous wave equation}) and (\ref{linearized kdv}),
Proposition \ref{mods behave} does not tell us that the solutions to (\ref{gun}) are controllable.  Nonetheless,
equation (\ref{gun}) is an inhomogeneous wave
equation where the inhomogeneity consists entirely of terms which are products
of left and right moving rapidly decaying functions.  From Wayne and Wright \cite{wayne.etal:02}
we have the following Lemma:
\begin{lemma}
\label{cross-term}
Suppose
$$
\partial_\tau u \pm \partial_\beta u = l(\beta+\tau,\epsilon^2 \tau) r(\beta-\tau,\epsilon^2 \tau),\quad u(X,0)=0.
$$
with $\| l(\cdot,T) \|_{H^s(2)} \le C$ and  $\| r(\cdot,T) \|_{H^s(2)} \le C$ for $T \in [0,T_0]$,
then
$$
\|u(\beta,\tau)\|_s \le C
$$
for all $\tau\in[0,T_0\ep^{-2}]$.  The constant $C$ is uniform in $\ep$.
\end{lemma}
Thus $W_3$ will remain $O(1)$.

\begin{remark}  
\label{more mods behave}
If we are in the situation in which 
Proposition \ref{mods behave} applies, we see that
least regular part in the driving term are $\partial_\beta^2 (UG)$ and $\partial_\beta^2(VF)$, which are
in $H^{\sigma-6}(2)$.  Thus, by this Lemma we have that $W_3 \in H^{\sigma-5}$ for 
all times of interest.
\end{remark}

At this time we have derived the modulation equations and shown the residual is small.  The only remaining order 
of business in this section is to determine how smooth the solutions to our modulation equations need to be in
order for $\resp$ to be appropriately regular.  This may seem to be a fairly tiresome task, but fortunately
the least regular terms in all of the sundry $E$ functions come from only one term---$E^u_{small}$!  This is because
$E^u_{small}$ contains many time derivatives.

We need to control $\res_a$ in $H^{s-1}$.  For this we need $\pt \res_u \in H^{s-1}$, which in turn
implies that we must have $\pt E^u_{small} \in H^{s-1}$.  Recalling equation (\ref{troublemaker}),
we see that this will require $\ptau^3 Z_2 \in H^{s-1}$, or rather (since $\Lep^{-1}$ saves a derivative)
$\ptau^3 W_2 \in H^{s-2}$.  
For this, we need $\pT^3 F$ and $\pT^3 G$ in $H^{s-2}$.  Given that $F$ solves (\ref{linearized kdv})
where $J^+$ contains the terms $\pbm^5 U$, $\pbm^3 \varphi^-$, one sees that $\pT^3 F$ will include the terms
$\pbm^9 F$, $\pbm^{11} U$ and $\pbm^9 \varphi^-$.  So $\|\pT^3 F\|_{s-2}$ is controlled by the $H^{s+9}$
norms of $U$ and $V$, and the $H^{s+7}$ norms of $\varphi^-$, $F$ and $G$.  The analogous result is
true for $\pT^3 G$.  
We also need $\ptau^3 W_3 \in H^{s-2}$.  Since $W_3$ solves (\ref{gun}), we require $W_3 \in H^{s+1}$.

In summary we have the following
Proposition.
\begin{proposition}
\label{residual}
Take $\pd$ as in (\ref{psi d}), with $U$, $V$, $F$, $G$, $P$
and $W_3$ solving their respective equations.  Let 
\begin{equation*}
\begin{split}
\pz = & L^{-1} \pd \\
\py = & \pz + \ep^4 \Delta_1 + \ep^6 \Delta_2\\
\pu = & \pa^{-1} \pt \pd\\
\Psi^a = & \pa^{-1} \pt^2 \pd
\end{split}
\end{equation*}
with $\Delta_1$ and $\Delta_2$ as in (\ref{delta 1}) and (\ref{delta 2}),
and form $\resp$ as in (\ref{res def 1}).  Then
\begin{equation*}
\begin{split}
\|\res_z\|_s \le & C \ep^{17/2}\\
\|\res_y\|_s \le & C \ep^{17/2}\\
\|\res_u\|_{s-1} \le & C \ep^{17/2}\\
\|\res_a\|_{s-1} \le & C \ep^{19/2}.
\end{split}
\end{equation*}
where $C$ is a constant which depends on
$\|U\|_{s+9}$,
$\|V\|_{s+9}$,
$\|P\|_{s+7}$,
$\|F\|_{s+7}$,
$\|G\|_{s+7}$ and 
$\|W_3\|_{s+1}$.  
The estimate (\ref{HOPE}) holds as long as these quantities remain bounded.  The constant
$C$ does not depend on $\ep$.
\end{proposition}

In light of Proposition \ref{mods behave} and Remark \ref{more mods behave}, we see that if
in we take the initial conditions for $U$ and $V$ to satisfy (\ref{ic kdv bound}) with $\sigma \ge s+7$,
that $\|U\|_{s+9}$,
$\|V\|_{s+9}$,
$\|P\|_{s+7}$,
$\|F\|_{s+7}$,
$\|G\|_{s+7}$ and 
$\|W_3\|_{s+1}$ are all $O(1)$ for $t \in [0,T_0 \ep^{-3}]$. And so we move on.

\section{The operator $K(x,y)$ Part II:  Estimates and Extensions}
\label{operator K 2}

In this section we will describe a few more estimates related to $K(x,y)$.  All
such estimates are either smoothing estimates or ones which show that
certain terms are small in the long wavelength setting.  

First, since $1 + \widehat{K_0}^2(k)$ goes to zero exponentially fast as
$|k| \to \infty$, the operator $1 + K_0^2$ is smoothing.  That is, for all $s\ge 0$
$$
\|(1+K_0^2) u\|_s \le C \|u\|_{L^2}.
$$

Also, commutators involving $K_0$ are smoothing.  We quote the following Lemma from 
\cite{schneider.etal:00}. 
\begin{lemma}
\label{commutator}
Let $r \ge 0$, $q > 1/2$, and $0 \le p \le q$.  Then there exists a $C > 0$  such that
$$
\| [f,K_0] g \|_r \le C \|f\|_{r+p} \|g\|_{q-p}.
$$
\end{lemma}
\begin{proof} See Lemma 3.12 on {\it p} 1498 of \cite{schneider.etal:00}.  \end{proof}

Schneider and Wayne show that $K_1(z,y)$ is a smoothing operator. 
\begin{lemma} 
\label{M_1 is bounded}
For $r \ge 0$, $q \ge 1/2$ and $0 \le p \le q$,
there is $C$ such that
\begin{equation}
\|K_1(z,y) u\|_r \le C \left(\|z\|_{r+p} + \|y\|_{r+p}\right) \|u\|_{q-p}.
\end{equation}
\end{lemma}
\begin{proof}
See Corollary 3.13 on {\it p} 1499 of \cite{schneider.etal:00}.  
\end{proof}

If we let $S_1(z,y) = K(z,y) - K_0$, we also have the following
estimates from \cite{schneider.etal:00}:
\begin{lemma}
\label{S_i facts} 
Fix $s \ge 4$. If the free surface is sufficiently smooth, then for $j=1,2$ we have:
\begin{itemize}
\item 
$$
\|S_j(z,y) u\|_s \le C \left(\|z\|^j_s + \|y\|^j_s\right) \|u\|_3,
$$
that is, $S_j$ is a smoothing operator,
\item
$$
\|\pa (S_j(z,y) u)\|_s \le C \left(\|z\|^j_s + \|y\|^j_s\right) \|u\|_3,
$$
that is, $\pa S_j$ is a smoothing operator,
\item
$$
\|[\pt,S_j] u\| \le C \left(\|z\|^j_s + \|y\|^j_s\right) \|u\|_3,
$$
that is, $[\pt,S_j]$ is a smoothing operator and this operator can be bounded
independently of $\pt u$ and
\item
$$
\|[\pt^2,S_j] u\| \le C \left(\|z\|^j_s + \|y\|^j_s\right) \left(\|u\|_4 + \|\pt u\|_4\right),
$$
that is, $[\pt^2,S_j]$ is a smoothing operator and this operator can be bounded
independently of $\pt^2 u$.
\end{itemize}
\end{lemma}
\begin{proof}
In \cite{schneider.etal:00}, see Lemmas 3.14, 3.15 and Corollary 3.16 on {\it pp} 1500, 1506 and 1507 respectively.
\end{proof}

We will also need the following propositions concerning the behavior of the remainder terms $S_1$ and $S_2$.
The first of these says that more or less the remainder $S_2$ is negligible for the sort of scalings we are considering.
That is to say, the term $E^y_{S_2}$ in the Section \ref{derivation} is very small.  
\begin{proposition}
\label{S_2 special}
Fix $s>5/2$.  Suppose $z=\ep^2 Z(\ep \al)$, $y=\ep^2 Y(\ep \al)$ and $f = \ep^2 F(\ep \al)$, with
$Z,Y,F \in H^{s+1}(2)$.  Moreover, assume $z-y = \ep^4 \Delta(\ep \al)$ with
$\Delta \in H^{s+1}(2)$.  
Then there exist $\ep_0$ such that for $\ep \in [0,\ep_0]$ 
there is a constant $C$ independent of $\ep$ such that:
$$
\| S_2(z,y) f \|_s \le C \ep^{17/2}.
$$
\end{proposition}

The second is a technical version of the mean value theorem as applied to the operator $S_1$.  

\begin{proposition}
\label{MVT} 
Suppose $z(\al,t)=\ep^2 Z(\ep (\al \pm t),\ep^3 t)$,  
$y(\al,t)=\ep^2 Y(\ep (\al \pm t),\ep^3 t)$, 
$u(\al,t)=\ep^2 U(\ep (\al \pm t),\ep^3 t)$ and   
$f(\al,t)=\ep^2 F(\ep (\al \pm t),\ep^3 t)$ with $Z,Y,U,F \in H^{s}(2)$ for 
$t \in [0,T_0 \ep^{-3}]$.  Also suppose
$R^z(\al,t)$, $R^y(\al,t)$ and $R^u(\al,t)\in H^s$ for the same time interval.  Then
\begin{align*}
\|  S_1(z(\cdot)+\ep^{11/2}R^z(\cdot),y(\cdot)+\ep^{11/2}R^y(\cdot)) f(\cdot) 
  - S_1(z(\cdot),y(\cdot)) f(\cdot) \|_s \le C \ep^{17/2}
\end{align*}
for $t\in[0,T_0 \ep^{-3}]$.
\end{proposition}

\begin{proof} {\bf for Proposition \ref{S_2 special}}:

First, notice $x(\al) = \int_0^\al Lz(a) da = \ep X(\ep \al)$.  
We know that $X$ is in $L^{\infty}$ by the following Lemma.
\begin{lemma}
\label{weight bound}
Suppose $f(\al) = \ep^2 F(\ep \al)$ with $F \in H^s(2)$.  Then for all $\al$
$$
\big | \int_0^{\al} f(a) da \big | \le C \ep \|F\|_{H^s(2)}
$$
\end{lemma}
\begin{proof} See Section \ref{proofs}.\end{proof}

Let $\Phi(\tilde{x},y)=(\Phi_1,\Phi_2)$ 
be the analytic map which takes $\Omega(t)$ 
to $$P^{-} = \left\{(\xi,\gamma) | \gamma \in [-1,0] \right\}.$$  That such a map exists
and is analytic is guaranteed by the Riemann mapping theorem.  
Let 
$$
h(\al) = \Phi_1(\tilde{x}(\al),y(\al))
$$ 
and $Q f = f \circ h$.  From \cite{schneider.etal:00}, we know that 
\begin{equation}
\label{explicit K}
K(x,y) f (\al) = Q \circ K_0 \circ Q^{-1} f (\al).
\end{equation}

We can derive a very useful implicit formula for $h^{-1}$ as follows.  The function
$\Phi^{-1}(\xi,\gamma)$ is analytic on $P^{-}$, thus it satisfies the Cauchy-Riemann
equations.  If we set
$$
\Phi^{-1}(\xi,\gamma) = \left(\xi + u_1(\xi,\gamma)\right) + i\left(\gamma -  v_1 (\xi,\gamma)\right)
$$
and notice that $\Phi^{-1}$ sends the bottom and top of $P^{-}$ to the bottom and top of $\Omega(t)$
respectively, we see that we have the following system
\begin{equation*}
\label{cauchy riemann}
\begin{split}
&\partial_\xi u_1 + \partial_\gamma v_1  = 0, \\
&\partial_\gamma u_1 - \partial_\xi v_1  = 0,\\
v_1(\xi,-1) &= 0, \quad
v_1(\xi,0)  = \eta(\xi)
\end{split}
\end{equation*}
where $\eta(\xi) = y(h^{-1}(\xi))$.  One can solve this system with the use of Fourier transforms
relatively simply.  If one does so, one finds that 
\begin{equation*}
\label{cauchy riemann solution}
u_1(\xi,0) = -\int_0^{\xi} \eta(\xi_1) d \xi_1 - M \eta (\xi)
\end{equation*}
where $M$ is the pseudo-differential operator given by
$$
\widehat{M \eta}(k) = \frac{k \cosh(k) - \sinh(k)}{ik \sinh(k)} \widehat{\eta}(k).
$$
Notice that to lowest order, $M$ is  $C \pa$.  

Now, notice that $h^{-1}(\xi) = \tilde{x}^{-1}(u(\xi,0))$ and so we have an 
implicit equation for $h^{-1}$.  
\begin{equation}
\label{form of h inverse}
\begin{split}
h^{-1}(\xi) = &\tilde{x}^{-1}\left(\xi -\int_0^{\xi} \eta(\xi_1) d \xi_1 - M \eta (\xi)\right)\\
       = &\tilde{x}^{-1}\left(\xi -\int_0^{\xi} y(h^{-1}(\xi_1)) d \xi_1 - M(y \circ h^{-1})  (\xi)\right).
\end{split}
\end{equation}
\begin{remark}
In Schneider and Wayne \cite{schneider.etal:00}, on {\it p} 1494, they make a minor
error in calculating this same function.  As a result, they claim the above representation
gives an explicit formula for $h^{-1}$.  Our correction here changes nothing about subsequent
steps in their proofs.  
\end{remark}

Since, $\tilde{x} = \al + \ep X(\ep \al)$ where $X$ is well-behaved, we can expect a similar
form for $\tilde{x}^{-1}$.  

\begin{lemma}
\label{inverse lemma}
Suppose $f(\al) = \al + g(\al)$ with $\|g\|_{C^{2}} \le 1/2$.  Then 
$$
f^{-1}(\xi) = \xi - g(\xi) + g(\xi) g^{\prime}(\xi) + E
$$
where $E = O(\|g\|^3_{W^{2,\infty}})$.  More specifically
$$
E \le C\left(\|g^{\prime}\|^2_{L^{\infty}} \|g\|_{L^{\infty}} 
           + \|g\|^2_{L^{\infty}} \|g^{\prime \prime}\|_{L^{\infty}} \right).
$$
In particular, notice that if $g(\al) = \ep G(\ep \al)$ this means $E = O(\ep^5)$.  
\end{lemma}
\begin{proof} See Section \ref{proofs}.\end{proof}

We apply this Lemma to $\tilde{x}$, and find that
\begin{align*}
\tilde{x}^{-1}(\xi) = \xi - \ep X(\ep \xi) + \ep^3 X(\ep \xi)\pb X(\ep \xi) + O(\ep^5).
\end{align*}
Combining this with (\ref{form of h inverse}) we can 
determine $h^{-1}(\xi)$ (and therefore $h$) in terms
of $x$ and $y$ to any order we wish.  To lowest order we see that 
\begin{align}
\label{snake}
h^{-1}(\xi) = \xi + O(\ep)
\end{align}
So now we have
\begin{equation*}
\begin{split}
h^{-1}(\xi)= &\xi -\int_0^{\xi} y(h^{-1}(\xi_1)) d \xi_1 - M(y \circ h^{-1})(\xi)\\
             &-\ep X\left(\ep \left(\xi -\int_0^{\xi} y(h^{-1}(\xi_1)) d \xi_1 - M(y \circ h^{-1})(\xi)\right)\right)\\
             &+ O(\ep^3)\\
           = &\xi -\int_0^{\xi} \ep^2 Y(\ep h^{-1}(\xi_1)) d \xi_1 - M(\ep^2 Y \circ \ep h^{-1})(\xi)\\
  &-\ep X\left(\ep \left(\xi -\int_0^{\xi} \ep^2 Y(\ep h^{-1}(\xi_1)) d \xi_1 
             - M(\ep^2 Y \circ \ep h^{-1})(\xi)\right)\right)\\
             &+ O(\ep^3).
\end{split}
\end{equation*}
If we insert (\ref{snake}) into the above 
and expand we have
\begin{equation*}
\begin{split}
h^{-1}(\xi) = \xi -\ep X(\ep \xi)  -\int_0^{\xi} \ep^2 Y(\ep\xi_1)) d \xi_1 - \ep^2 M (Y (\ep\cdot))(\xi)
+ O(\ep^3).
\end{split}
\end{equation*}

One can continue in this manner and determine the next order terms in the expansion
of $h^{-1}$.  
If we let
$$
\ep G_1(\ep \al) =  -\ep X(\ep \al) -\int_0^{\al} \ep^2 Y(\ep\al_1)) d \al_1 - \ep^2 M (Y (\ep\cdot))(\al).
$$
the expansion is
$$
h^{-1}(\xi) = \xi + \ep G_1(\ep \xi) + \ep^3 B_1(\ep \xi) + O(\ep^5)
$$
where
\begin{align*}
B_1(\xi) = &\int_0^{\xi} \ep G_1(\ep a) \pb \ep^3 Y (\ep a) da\\
           & + M( \ep^3 G_1(\ep \cdot) \pb Y(\ep \cdot)) (\xi) + 
           \ep^3 G_1(\ep \xi) \pb X (\ep \xi).
\end{align*}

Notice that since $M$ is $C \pa$ to lowest order, $- \ep^2 M (Y (\ep\cdot))(\al)$ is $O(\ep^3)$.
Moreover,
by hypothesis, we have $\ep^2 Z(\ep \al) - \ep^2 Y(\ep \al) = \ep^4 \Delta(\ep \al)$.  Thus
\begin{align*}
 - \ep X(\ep \al)-\int_0^{\al} \ep^2 Y(\ep a) da 
=&-\int_0^{\al} \left(\ep^2 \pb X(\ep a) + \ep^2 Y(\ep a)\right) da\\
=&-\int_0^{\al} \left(\ep^2 L(Z(\ep \cdot))(a) + \ep^2 Y(\ep a)\right) da\\
=&\int_0^{\al} \left(\ep^2 Z(\ep a) - \ep^2 Y(\ep a) \right) da + O(\ep^3) \\
=&\int_0^{\al} \left(\ep^4 \Delta(\ep a) \right) da + O(\ep^3) \\
=&O(\ep^3).  
\end{align*}
That is, $\ep G_1$ is really $O(\ep^3)$!  This cancelation is the crucial step in this
proof.
Since $\ep G_1$ appears in each term in $B_1$, we have shown
$$
h^{-1}(\xi) = \xi + \ep^3 G(\ep \xi) + O(\ep^5)
$$
with $\ep^3 G = \ep G_1$.  
We appeal to Lemma \ref{inverse lemma} again, and we have
$$
h(\al) = \al - \ep^3 G(\ep \al) + O(\ep^5).  
$$

Now that we have particularly good estimates on $h$ and $h^{-1}$, we can begin our
discussion of $K$ in earnest. For notational simplicity, we will let
\begin{equation*}
\begin{split}
h(\al) = \al + g_1(\al)\\
h^{-1}(\xi) = \xi + g_2(\xi)
\end{split}
\end{equation*}
If we let $$\tilde{f} = Q^{-1} f$$  
we can make the following formal approximation
using Taylor's theorem,
\begin{equation*}
\begin{split}
K(x,y) f (\al) = & Q \circ K_0 \tilde{f} (\al)\\
               = & K_0 \tilde{f} (h(\al))\\
               = & K_0 \tilde{f} (\al + g_1(\al))\\
               = & K_0 \tilde{f} (\al) + g_1(\al) K_0 \pa \tilde{f} (\al) + h.o.t..
\end{split}
\end{equation*}
Also by Taylor's theorem,
$$
\tilde{f}(\al)= f(\al) + g_2(\al) \pa f (\al) + h.o.t..
$$
Putting these together we have
\begin{align*}
K(x,y) f(\al)=& K_0 f (\al) + g_1(\al) K_0 f^\prime (\al) + K_0(g_2 f^\prime) (\al) + h.o.t..
\end{align*}
So let
\begin{equation*}
\begin{split}
E_1 f =& K_0 f (\al) + g_1(\al) K_0 f^\prime (\al) + K_0(g_2 f^\prime) (\al)\\
E_2 f =& K(x,y) f - E_1 f.
\end{split}
\end{equation*}
We prove Proposition \ref{S_2 special} if we can prove
\begin{itemize}
\item 
$\| E_1 f -  K_0 f - K_1(x,y) f \|_s \le  C\ep^{17/2}$, and
\item
$\|E_2 f\|_s \le C \ep^{17/2}$.
\end{itemize}

Let us deal with $E_2 f$ first.  We can rewrite $E_2 f$ as:
$$
E_2 f  = E_2^1 f + E_2^2 f +E_2^3 f
$$
with
\begin{equation*}
\begin{split}
E_2^1 f  = &K(x,y) f - K_0 \tilde{f} - g_1 K_0 \pa \tilde{f} \\
E_2^2 f  = &K_0 \tilde{f} - K_0 f - K_0(g_2 \pa f)\\ 
E_2^3 f  = &g_1 K_0 \pa \tilde{f} - g_1 K_0 \pa f \\
\end{split}
\end{equation*}

As our approximation for $K$ was determined by an application of Taylor's theorem,
we need to prove a Lemma which shows that this formal step 
can be made rigorous, at least for functions in the weighted Sobolev spaces.

\begin{lemma}
\label{the weight}
Suppose $F \in H^{s}(n)$, $s>1/2$, $n>1/2$. Then for all $C_0>0$ there exists $\ep_0$ such
that for $\ep \in [0,\ep_0]$ there is a constant $C$ independent of $\ep$ such that:
$$
\left ( \int_{|\al| > C_0 \ep^{-3}} |F(\ep \al)|^2 d\al \right)^{1/2} \le C \ep^{2n - 3/2}.
$$ 
Moreover, for $1 \le j \le s$ we have:
$$
\left ( \int_{|\al| > C_0 \ep^{-3}} |\pa^j F(\ep \al)|^2 d\al \right)^{1/2} \le C \ep^{2n-3/2+j}.
$$ 
\end{lemma}
\begin{proof} See Section \ref{proofs}.\end{proof}

\begin{remark}
If instead we are considering 
$$
\left ( \int_{|\al| > C_0 \ep^{-3}} |F(\ep (\al \pm t),\ep^3 t)|^2 d\al \right)^{1/2} 
$$
with $F(\cdot,T) \in H^s(n)$ for $T \in [0,T_0]$, we can maintain the same bound
as above by taking $C_0 \ge 2 T_0$.  
\end{remark}

We can use the above Lemma to prove a version of Taylor's theorem.
\begin{lemma}
\label{mean value theorem}
Suppose $F \in H^s(2)$, $s>5/2$ and $g \in L^\infty$.  Then
$$
\| F(\ep \cdot + \ep^2 g(\ep \cdot)) - F(\ep \cdot)\|_{L^2} \le C \ep^{3/2}.
$$
\end{lemma}

\begin{proof}  
By Lemma \ref{the weight} we have
\begin{equation*}
\begin{split}
   &\| F(\ep \cdot + \ep^2 g(\ep \cdot)) - F(\ep \cdot)\|^2 \\
=   &\int_{|\al| \le \ep^{-3}} |F(\ep \al + \ep^2 g(\ep \al)) - F(\ep \al)|^2 d\al \\
   +&\int_{|\al| \ge \ep^{-3}} |F(\ep \al + \ep^2 g(\ep \al)) - F(\ep \al)|^2 d\al \\
\le &\int_{|\al| \le \ep^{-3}} |F(\ep \al + \ep^2 g(\ep \al)) - F(\ep \al)|^2 d\al  + C\ep^{5}.
\end{split}
\end{equation*}
Now, we add and subtract $\ep^2 g(\ep \al) F^{\prime}(\ep \al)$ in the
remaining integral,
\begin{equation*}
\begin{split}
   &\int_{|\al| \le \ep^{-3}} |F(\ep \al + \ep^2 g(\ep \al)) - F(\ep \al)|^2 d\al\\
\le&\int_{|\al| \le \ep^{-3}} |F(\ep \al + \ep^2 g(\ep \al)) - F(\ep \al)- \ep^2 g(\ep \al) F^{\prime}(\ep \al)|^2 d\al\\
   &+\int_{|\al| \le \ep^{-3}} |\ep^2 g(\ep \al) F^{\prime}(\ep \al)|^2 d\al\\
\le&\int_{|\al| \le \ep^{-3} } 
|F(\ep \al + \ep^2 g(\ep \al)) - F(\ep \al) - \ep^2 g(\ep \al) F^{\prime}(\ep \al)|^2 d\al \\
   &+ \ep^{3} \|g(\cdot)\|^2_{L^{\infty}} \|F^{\prime}(\cdot)\|_{L^2}^2.
\end{split}
\end{equation*}
We naively bound the above integral and apply the mean value theorem.  That is,
\begin{equation*}
\begin{split}
&\int_{|\al| \le \ep^{-3}} 
|F(\ep \al + \ep^2 g(\ep \al)) - F(\ep \al) - \ep^2 g(\ep \al) F^{\prime}(\ep \al)|^2 d\al\\
\le &C \ep^{-3} \sup_{|\al| \le \ep^{-3}} 
|F(\ep \al + \ep^2 g(\ep \al)) - F(\ep \al) - \ep^2 g(\ep \al) F^{\prime}(\ep \al)|^2 \\
\le &C \ep^{-3} \sup_{|\al| \le \ep^{-3}} |\ep^4 g^2(\ep \al) F^{\prime \prime}(\ep \al^*)|^2\\
\le &C \ep^{5} \|g\|^4_{L^\infty} \|F^{\prime \prime}\|^2_{L^\infty}.
\end{split}
\end{equation*}
With this, we have proven the Lemma.  
\end{proof}

\begin{remark}
With this general technique we are also able to show that
$$
\| \pa^j(F(\ep \cdot + \ep^2 g(\ep \cdot)) - F(\ep \cdot))\|_{L^2} \le C \ep^{3/2+j},
$$
$$
\| \pa^j\left(F(\ep \cdot + \ep^2 g(\ep \cdot)) - F(\ep \cdot) - \ep^2 g(\ep \cdot) F^\prime(\cdot)\right) \|_{L^2} \le C \ep^{7/2+j},
$$
\begin{align*}
\|& \pa^j\left(F(\ep \cdot + \ep^2 g(\ep \cdot)) \right)\\
  & -\pa^j\left(  F(\ep \cdot) + \ep^2 g(\ep \cdot) F^\prime(\cdot) +1/2 \ep^4 g^2(\ep\cdot) F^{\prime\prime}(\ep \cdot) \right)\|_{L^2} 
\le C \ep^{11/2+j}
\end{align*}
and so on.  
\end{remark}

Now we will be able to control $E^1_2$.  We can control the other functions in precisely the same
fashion.  Since $f$ is long wavelength and of rapid decay, so is $\tilde{f}$.  And thus we 
can use Lemma \ref{mean value theorem}.  
In what follows, $\ep^3 \tilde{F} (\ep \al) = K_0 \tilde{f}(\al)$. (The
extra $\ep$ comes from the long wave approximation of $K_0$.)  
\begin{equation*}
\begin{split}
 & \| E^1_2 f (\cdot) \|_s\\
 =&\|K_0 \tilde{f} (\cdot + g_1(\cdot))  
   - K_0 \tilde{f}(\cdot) - g_1(\cdot) K_0 \pa \tilde{f}(\cdot)\|_s\\
\le& C \ep^3 \| \tilde{F}(\ep \cdot + \ep^4 G(\ep \cdot)) 
  - \tilde{F}(\ep \cdot) - \ep^4 G(\ep \cdot) \pb \tilde{F}(\ep \cdot)\|_s\\
\le&C\ep^{17/2} 
\end{split}
\end{equation*}

Now we turn our attention to 
$E_1 f -  K_0 f - K_1(x,y) f$.  A routine calculation shows that this is equal to
\begin{equation*}
\begin{split}
 &-\ep^3 [G(\ep \al),K_0] \ep^3 \pb F(\ep \al) - K_1(\ep X,\ep^2 Y)\ep^2 F(\ep \al)\\
 &+(g_1+\ep^3 G) K_0 \ep^3 \pb F(\ep \al) + K_0 \left( (g_2 - \ep^3 G) \ep^3\pb F(\ep \cdot) \right)(\al)
\end{split}
\end{equation*}
Now, $g_1 + \ep^3 G$ and $g_2 - \ep^3 G$ are $O(\ep^5)$, so the second line above can be
bounded by $C\ep^{17/2}$, if we make use of $K_0$'s long wavelength
approximation.  Moreover, we claim that the first line is identically zero.
Let
\begin{equation*}
\begin{split}
b(\al) = &-\int_0^{\al} \ep^2 Y(\ep\al_1) d \al_1 - \ep^2 M (Y (\ep\cdot))(\al)\\
       = &-\int_0^{\al}  y(\al_1) d \al_1 - M y (\al).
\end{split}
\end{equation*}
Thus $\ep^3 G(\ep \al) = b(\al) - x(\al)$.  So
$$
\ep^3 [G(\ep \al),K_0] \ep^3 \pb F(\ep \al) = -[x(\al),K_0] \pa f(\al) + [b(\al),K_0] \pa f(\al).
$$
Taking the Fourier transform of the second term we have
\begin{align*}
=&\mathfrak{F} \left( [b(\al),K_0] \pa f(\al) \right) (k)\\
=&\int \left(\widehat{K_0}(l)-\widehat{K_0}(k)\right) \widehat{b}(l-k) il \widehat{f}(l) dl\\
=&\int \frac{\widehat{K_0}(l)-\widehat{K_0}(k)}{i(l-k)} i(l-k) \widehat{b}(l-k) il \widehat{f}(l) dl\\
=&\int \frac{\widehat{K_0}(l)-\widehat{K_0}(k)}{i(l-k)} \widehat{\pa b}(l-k) il \widehat{f}(l) dl
\end{align*}
Now notice that
$
\pa b = L y,
$
so the above becomes
\begin{align*}
 &\int \frac{\widehat{K_0}(l)-\widehat{K_0}(k)}{i(l-k)} \widehat{L y}(l-k) il \widehat{f}(l) dl\\
=&\int \frac{\widehat{K_0}(l)-\widehat{K_0}(k)}{\widehat{K_0}(l-k)} \widehat{y}(l-k) il \widehat{f}(l) dl\\
=&\int (1 + \widehat{K_0}(k) \widehat{K_0}(l)) \widehat{y}(l-k) il \widehat{f}(l) dl\\
=&\mathfrak{F}\left( (y + K_0 y K_0) \pa f\right)(k).
\end{align*}
where we have used the trigonometric identity (\ref{trig id}) from Section \ref{derivation}.

That is
\begin{align*}
  &\ep^3 [G(\ep \al),K_0] \ep^3 \pb F(\ep \al) \\
= &-[x(\al),K_0] \pa f(\al) + \left( (y + K_0 y K_0) \pa f\right)(k)\\
= &-K_1(x,y) f
\end{align*}
and so we are done with the proof of Proposition \ref{S_2 special}.  

\end{proof}

\begin{proof} {\bf for Proposition \ref{MVT}}:

Let $h$ be as in the above subsection and 
$h_2$ be the analogous function 
for the configuration $(z_2,y_2)=(z+\ep^{11/2} R^z,y+\ep^{11/2} R^y)$.  We also define
the function $\tilde{x}_2 = \al + x_2$ by $\pa x_2 = L z_2$.  Unlike in the previous Lemma,
here the time dependence of the functions is important.  And so we determine $x$ and $x_2$ by
integrating in both space and time.  That is,
\begin{equation*}
\begin{split}
\tilde{x}(\al,t)=&\alpha + \epsilon \chi\left(\ep \alpha,\ep t\right)\\
\tilde{x}_2(\al,t)=&\alpha + \epsilon \chi\left(\ep \alpha,\ep t\right) + \epsilon^{5/2} E(t) + \ep^{11/2} \rho(\al,t).
\end{split}
\end{equation*}
with
\begin{equation*}
\begin{split}
\ep \chi(\ep \al, \ep t) = & \left( \int_0^t u\left(0,s\right) ds + \int_0^\alpha L z\left(w,t\right) dw \right)\\
\ep^{5/2}E(t) = &  \epsilon^{11/2} \int_0^t   R^u\left(0,s\right) ds\\ 
\ep^{11/2} \rho = & \int_0^\alpha   L R^z\left(w,t\right) dw 
\end{split}
\end{equation*}

The functions satisfy the following estimates for all $t \in [0,T_0 \ep^{-3}]$.
\begin{equation*}
\begin{split}
|\ep \chi(\ep \al, \ep t)| &\le  C \ep\\
|\ep^{5/2}E(t) | &\le  C \ep^{5/2}\\
|\ep^{11/2}\rho(\al,t)| &\le  C\sqrt{|\al|} \|R^z\|_{H^s}
\end{split}
\end{equation*}

The first estimate follows from similar estimates in the previous Lemma, the second is the 
naive bound and the final follows from following simple fact:
\begin{fact}
If $p(\al) = \int_0^{\al} r(a) da$, where $r \in L^2$, then
$$
|p(\al)| \le \sqrt{|\al|} \|r\|_{L^2}.
$$
\end{fact}
In what follows we will make strong use of the fact that $E(t)$ does not depend on $\al$.

Using the same techniques as were used in proving Lemma \ref{inverse lemma}, one
can show that:
$$
\tilde{x}^{-1}_2 (\xi,t) = \tilde{x}^{-1}(\xi,t) - \ep^{5/2} E(t) + \ep^{11/2} \rho_2(\xi)
$$
where $|\rho_2(\al)| \le \sqrt{|\al|} \|R^z\|_{L^2}$.  This sort of 
estimate carries over to the functions $h$.  That is
\begin{equation*}
\begin{split}
h_2^{-1}(\xi,t) = &h^{-1}(\xi,t) - \ep^{5/2} E(t) + \ep^{11/2} \rho_3(\xi,t)\\
h_2(\al,t)=  &h(\al,t) + \ep^{5/2} E(t) + \ep^{11/2} \rho_4(\al,t)
\end{split}
\end{equation*}
where $|\rho_3(\al,t)|,|\rho_4(\al,t)| \le \sqrt{|\al|} (\|R^z\|_{L^2}+\|R^y\|_{L^2})$
over the long time scale.

Now 
define $Q_2 f = f \circ h_2$.  So
\begin{equation*}
\begin{split}
 &S_1(z + \ep^{11/2} R^z,y + \ep^{11/2}) f - S_1(z,y) f\\
=&Q_2 \circ K_0 \circ Q_2^{-1} f - Q \circ K_0 \circ Q^{-1} f\\
=&Q \circ \left(Q^{-1} \circ Q_2 \circ K_0 \circ Q_2^{-1} \circ Q - K_0 \right) \circ Q^{-1} f
\end{split}
\end{equation*}
Since $Q$ and its inverse are bounded operators from $H^s$ to $H^s$, we need only prove the estimate for the
operator
$$
\tilde{Q}\circ K_0 \circ \tilde{Q}^{-1} - K_0 
$$
where $\tilde{Q} = Q^{-1} \circ Q_2$.  
Notice that $\tilde{Q} \circ K_0 \circ \tilde{Q}^{-1}$ is the Hilbert operator $K$ for a domain  
with the ``$h$'' function given by $\tilde{h}(\al)=h_2(h^{-1}(\al))$.  Moreover, from
the above calculations for $h$ and $h_2$ we have,
$$
h_2(h^{-1}(\al)) = \al + \ep^{5/2} E(t) + \ep^{11/2} \rho_5(\al,t)
$$
with $\rho_5$ satisfying the same type of estimates as $\rho_4$.  

At this point we can make an appeal to Lemma 3.14 on {\it p} 1500 of \cite{schneider.etal:00}.  
In this lemma they prove that
$\| S_1(z,y) f \|_s \le C \left(\|z\|_s + \|y\|_s\right) \| f \|_3$.
In the course of their proof, they show that if
$h(\al) = \al + g(\al)$ then
\begin{equation*}
\|Q\circ K_0 \circ Q^{-1} f(\cdot) - K_0 f(\cdot)\|_s \le C \|\pa g\|_{s-1} \|\pa f\|_2.
\end{equation*}
(See the inequalities in Cases I-IV on {\it pp} 1501-1506.)  
So if we set
$\tilde{g} = \ep^{5/2} E(t) + \ep^{11/2} \rho_5(\al,t)$, we see that taking a spatial
derivative leaves us with 
$$
\pa \tilde{g} =   O(\ep^{11/2}).
$$
And so, if we keep in mind that $f(\al,t) = \ep^2 F(\ep(\al\pm t,\ep^3 t)$,
\begin{equation*}
\begin{split}
\|\tilde{Q}\circ K_0 \circ \tilde{Q}^{-1} f(\cdot) - K_0 f(\cdot)\|_s \le &C \|\pa \tilde{g}\|_{s-1} \|\pa f\|_2\\
      \le &C \ep^{17/2}
\end{split}
\end{equation*}
This completes the proof of Proposition \ref{MVT}.
\end{proof}

\section{The error estimates}
\label{error estimates}
In this section we prove that the approximation is rigorous.  That is
we will prove Theorem \ref{main result water wave}.  
We will be working with the three
and four dimensional formulations of the water wave problem (equations (\ref{three-d water wave})
and (\ref{four-d water wave})).  From \cite{schneider.etal:00}, we know that for initial data 
of the type we are considering, solutions to these equations exist over the long times
we are considering.  If $(z,y,u)$ is a solution to (\ref{three-d water wave}),
let:  
\begin{equation}
\label{error functions}
\begin{split}
z(\al,t) = &  {\Psi}^z( \al,  t) + \ep^{11/2} R^z(\al,t)\\
y(\al,t) = & {\Psi}^y( \al,  t) + \ep^{11/2} R^y(\al,t)\\
u(\al,t) = &  {\Psi}^u(\al,  t) + \ep^{11/2} R^u(\al,t)
\end{split}
\end{equation}
with the functions $\Psi$ defined as above.
We call $R^z$, $R^y$ and $R^u$ ``error'' functions and we
denote $\bar{R} = (R^z,R^y,R^u)$.  
Our goal will be to show that $\bar{R}$ 
remain $O(1)$ in $\sob^s = H^s \times H^s \times H^{s-1/2}$ over the long time
scale.  
If we can do this, then we will have proven the main theorem.  The first step will be
to determine the equations which these functions satisfy.  Loosely, we want to be able
to write for each of the error functions an evolution equation of the form
$$
\pt R = \textrm{quasilinear} + \textrm{small and smooth}.
$$
We will at times go to great lengths to achieve this!

Clearly,
\begin{equation*}
\label{error evolution z}
\begin{split}
\pt R^z &= K_0 R^u. 
\end{split}
\end{equation*}
Finding the equations for $R^y$ and $R^u$ is a bit more complex.  First we focus
on $R^y$.  Substituting from (\ref{error functions}) into $\pt y = K(z,y) u$, we have
\begin{equation*}
\begin{split}
 &\pt R^y\\ 
=& \ep^{-11/2} \left( K(\pz + \ep^{11/2} R^z,\py + \ep^{11/2} R^y) 
                         \left(\putt + \ep^{11/2} R^u\right) - \pt \py \right)\\
        =&\phantom{+} K(\pz + \ep^{11/2} R^z,\py + \ep^{11/2} R^y) R^u \\
                   &+\ep^{-11/2}\left( K(\pz+\ep^{11/2}R^z,\py+\ep^{11/2}R^y) 
                                      \putt - \pt \py \right) \\
        =&\phantom{+} K_0 R^u + M_1(\pz) \pa R^u  - (\py + K_0(\Psi^y K_0)) \pa R^u + N^y       
\end{split}
\end{equation*}
where
\begin{equation*}
\begin{split}
&N^y\\ =&  \ep^{-11/2} \textrm{Res}_y \\
     &+  \ep^{-11/2} \left( \left(S_1(\pz+\ep^{11/2}R^z,\py+\ep^{11/2}R^y) -
                        S_1(\pz,\py) \right)\putt \right)\\ 
     & +\left(K(\pz+\ep^{11/2}R^z,\py+\ep^{11/2}R^y) - K_0 - K_1(\pz,\py) \right) R^u.
\end{split}
\end{equation*}
We claim that $N^y$ is ``small''.  That is, we have:
\begin{lemma}
For all $C_R>0$, there exists $\epsilon_0$ such that for all $\ep \in (0,\ep_0)$ 
and $t$ such that $\sup_{0 \le t^\prime \le t} \|\bar{R}(\cdot,t^\prime)\|_{\sob^s} \le C_R$
we have:
$$
\|N^y\|_s \le C\left( \ep^3  +  \ep^3 \|\bar{R}\|_{\sob^s} + \ep^{11/2} \|\bar{R}\|^2_{\sob^s} \right).
$$
\end{lemma}
\begin{proof}

First we remark that the approximating functions $\Psi$, and their derivatives are
all bounded over the long time scales.  Thus, we will not be keeping track of the
dependence of the norm of $N^y$ on the norms of these functions.  
By Proposition \ref{residual}, we know that 
$\|\ep^{-11/2} \textrm{Res}_y\|_s \le C \ep^3$.  

We can bound 
$$
\ep^{-11/2}\left(S_1(\pz+\ep^{11/2}R^z,\py+\ep^{11/2}R^y) - S_1(\pz,\py) \right)\putt 
$$
by Lemma \ref{MVT}.  

Finally,       
\begin{align*}
&\|\left(K(\pz+\ep^{11/2}R^z,\py+\ep^{11/2}R^y)-K_0-K_1(\pz,\py)\right) R^u\|_s\\
\le &C \left(\ep^3 \|\bar{R}\|_{\sob^s} + \ep^{11/2} \|\bar{R}\|^2_{\sob^s} \right)
\end{align*}
by the estimates on $K$ and its expansions which we saw in Sections \ref{operator K 1} and 
\ref{operator K 2} (in particular
Lemma \ref{M_1 is bounded} and
Lemma \ref{S_i facts}).
\end{proof}

Now we discuss $R^u$.  We know that 
\begin{equation}
\label{ww 1 1}
\pt u (1 + Lz) + \pa y (1 + \pt^2 y) = 0.
\end{equation}
We would like an evolution type equation for $R^u$.  Notice that since $\pt y = K(z,y) u$, there
is a ``hidden'' $\pt u$ in the term $\pt^2 y$.  Recall that the commutator
$[\pt,S_1(z,y)] u$ can be bounded independently of $\pt u$ (see Lemma \ref{S_i facts} 
in Section \ref{operator K 2}). 
Therefore, we can rewrite the above as
$$
(1 + Lz + \pa y K(z,y)) \pt u  + \pa y (1 + [\pt,S_1(z,y)] u) = 0.
$$
Substituting in for $u$ from (\ref{error functions}) the above becomes
\begin{equation*}
\begin{split}
0=&\left( 1 + Lz + \pa y K(z,y) \right) \pt \ep^{11/2} R^u \\
  &         + \left( 1 + Lz + \pa y K(z,y) \right) \pt \ep^{2} \putt  \\
  &         + \pa y \left(1 + [\pt,S_1(z,y)] \putt \right)\\
  &         + \pa y \left([\pt,S_1(z,y)] \ep^{11/2} R^u \right)  
\end{split}
\end{equation*}
or rather
\begin{equation*}
\begin{split}
0=&\left( 1 + Lz + \pa y K(z,y) \right) \pt \ep^{11/2} R^u \\
  &         + \left( 1 + Lz \right) \pt \ep^{2} \putt  \\
  &         + \pa y \left(1 + \pt (K(z,y) \putt) \right)\\
  &         + \pa y \left([\pt,S_1(z,y)] \ep^{11/2} R^u \right)  .
\end{split}
\end{equation*}
We rearrange this a bit, and break up $y$ and $z$.
\begin{equation*}
\label{come back here}
\begin{split}
0                  =&\left( 1 + Lz + \pa y K(z,y) \right) \pt R^u + \pa R^y\\
                    &+ L R^z \pt \putt + \pa  R^y \pt (K(z,y) \putt) + \pa y [\pt,S_1(z,y)]  R^u   \\
                    &+\ep^{-11/2}\left(\left( 1 + L \pz \right) \pt \ep^{2} \putt  
                     + \pa \py \left(1 + \pt (K(z,y) \putt) \right)\right).
\end{split}
\end{equation*}

The operator 
$$
A(z,y) =  (1 + Lz + \pa y K(z,y))
$$
is invertible since $K(z,y)$ is a bounded operator on $H^s$, 
provided $z$ and $y$ are small (which they are).  
Moreover, we can approximate $A^{-1}$ {\it via} the Neumann series.  Thus
the above equation can be rewritten as
\begin{equation*}
\pt R^u = -(1 - \ep^2 W_1) \pa R^y + N^u
\end{equation*}
where $N^u= N^u_1 + N^u_2$ and 
\begin{equation*}
\begin{split}
N^u_{1} =& 
-\ep^{-11/2} A^{-1} \left(\left( 1 + L \pz \right) \pt \ep^{2} \putt  
                    + \pa \py \left(1 + \pt (K(z,y) \putt) \right)\right)\\
N^u_2 = &-A^{-1}
  \left(L R^z \pt \putt + \pa  R^y \pt (K(z,y) \putt) + \pa y [\pt,S_1(z,y)]  R^u \right)  \\
      &+\left(-A^{-1} + (1 - \ep^2 W_1) \right) \pa R^y.
\end{split}
\end{equation*}

\begin{lemma}
For all $C_R>0$, there exists $\epsilon_0$ such that for all $\ep \in (0,\ep_0)$ 
and $t$ such that $\sup_{0 \le t^\prime \le t} \|\bar{R}(\cdot,t^\prime)\|_{\sob^s} \le C_R$
we have
$$
\|N^u\|_{s-1} 
\le C\left( \ep^3  +  \ep^3 \|\bar{R}\|_{\sob^s} + \ep^{11/2} \|\bar{R}\|^2_{\sob^s} \right).
$$
\end{lemma}
\begin{proof}
First we point out this estimate is in $H^{s-1}$.  The loss of regularity here is easily
seen.  Both $L R^z$ and $\pa R^y$ explicitly appear in $N^u_2$, and
are not smoothed by any operators.  Thus, losing this derivative is unavoidable.  In fact,
it is easy to see that the above estimates holds for $N^u_2$ by noting that $A^{-1}$, $K$
and $[\pt,S_1]$ are bounded operators.  

Bounding $N^u_1$ is also easily done once we recognize that this term is almost
exactly $\ep^{-11/2}\textrm{Res}_u$.  We have
\begin{equation*}
\begin{split}
 &\ep^{-11/2} A^{-1} \left(\left( 1 + L \pz \right) \pt \ep^{2} \putt  
                    + \pa \py \left(1 + \pt (K(z,y) \putt) \right)\right)\\
=&\ep^{-11/2}A^{-1} \left(\left( 1 + L \pz \right) \pt \ep^{2} \putt  
                    + \pa \py \left(1 + \pt (K(\pz,\py) \putt) \right)\right)\\& + N^u_3\\
=&\ep^{-11/2}\left(\pt \ep^{2} \putt + \pa \py \frac{1 +\pt (K(\pz,\py) \putt)}
                                                           {1 + L\pz}\right) + N^u_3 + N^u_4\\
=& \ep^{-11/2} \left(\pt \ep^{2} \putt + \pa \py \frac{1 +\pt^2 (\py)}{1 + L\pz}\right) 
+ N^u_3 + N^u_4+ N^u_5\\
=& \ep^{-11/2}\textrm{Res}_u+ N^u_3 + N^u_4+ N^u_5
\end{split}
\end{equation*}
where
\begin{equation*}
\begin{split}
N^u_3 = &\ep^{-11/2}A^{-1}\left(\pa \py \pt \left(K(z,y)\putt - K(\pz,\py)\putt\right)\right)\\
N^u_4 = &\ep^{-11/2}(A^{-1}-(1+L\pz)^{-1}) \\
        &\times \left(\left( 1 + L \pz \right) \pt \ep^{2} \putt  
                        + \pa \py \left(1 + \pt (K(\pz,\py) \putt) \right)\right)\\
N^u_5 =&\ep^{-11/2} \pa \py \frac{\pt (\textrm{Res}_y)}{1 + L\pz}.
\end{split}
\end{equation*}
We bound $N^u_3$ using mean value theorem arguments entirely analogous to those used
when bounding $N^y$.  To bound $N^y_4$, one observes that 
$$
\left( 1 + L \pz \right) \pt \ep^{2} \putt  
                        + \pa \py \left(1 + \pt (K(\pz,\py) \putt) \right)
$$
is very nearly $\textrm{Res}_u$ and is thus $O(\ep^{17/2})$.  
$N^u_5$ is clearly small, as it contains $\pt \textrm{Res}_y$.  
This completes the proof.
\end{proof}

We need to make analogous calculations for the four dimensional system.  Let
$$
a(\al,t) = \ep^3 \patt(\ep \al, \ep t) + \ep^{11/2} R^a(\al,t)
$$
and $\bar{R}_e = (R^z,R^y,R^u,R^a)$.  This extended set of error functions
lives in $\mathfrak{H}^s_e = H^s \times H^s \times H^{s-1/2} \times H^{s-1}$.

It is easy to see that
$$
\pt R^u = R^a
$$
but more difficult to determine the evolution of $R^a$.  
We begin by taking a 
time derivative of (\ref{ww 1 1}).
\begin{equation}
\label{time deriv ww}
\pt^2 u (1 + Lz) + \pt u \pa u + \pa \pt y (1 + \pt^2 y) +\pa y \pt^3 y = 0.
\end{equation}
Letting,
\begin{equation*}
\begin{split}
I =& \pt^2 u (1+Lz) + \pa y \pt^3 y\\
II=& \pa \pt y (1 + \pt^2 y)+\pt u \pa u 
\end{split}
\end{equation*}
(\ref{time deriv ww}) is $I + II = 0$.  

Manipulations very similar to those carried out in
determing $\pt R^u$, show that 
\begin{equation*}
\begin{split}
I =& A(z,y) \pt \ep^{11/2} R^a + \pa y [\pt^2, S_1(z,y)] \ep^{11/2} R^u\\
   & + (1 + Lz) \pt \ep^3 \patt + \pa y \pt^2 (K(z,y) \putt).
\end{split}
\end{equation*}

For $II$, we have
\begin{equation*}
\begin{split}
II = & (1 + \pt^2 y) \pa \pt \ep^{11/2} R^y + \pt u \pa \ep^{11/2} R^u \\
     & +(1 + \pt^2 y) \pa \pt \py + \pt u \pa \putt\\
   = & \ep^{11/2}\left( (1+\pt^2 y)\pa (K_0 R^u + K_1(\pz,\py)R^u + N^y)  + \pt u \pa R^u\right)\\  
     & +(1 + \pt^2 y) \pa \pt \py + \pt u \pa \putt\\
   = & \ep^{11/2} \left( (1 + \pt^2 y - \pt u K_0 \cdot) K_0 \pa R^u + \pa K_1 (\pz,\py) R^u \right)\\
     & +(1 + \pt^2 y) \pa \pt \py + \pt u \pa \putt + \ep^{11/2} B_{II}
\end{split}
\end{equation*}
where
$$
B_{II} = \pt u (1 + K_0^2) \pa R^u + \pt^2 y \pa (K_1(\pz,\py) R^u + N^y) + \pa N^y.
$$
Noting that $\pt^2 y = K_0 a + [\pt, S_1] u + S_1 a$, we see that
that $B_{II}$ is smooth in the error functions, and is $O(\ep^3)$.    

Adding $I$ and $II$ gives
\begin{equation}
\label{ishmael}
\begin{split}
0 =&A(z,y) \pt R^a +  (1 + \pt^2 y - \pt u K_0 \cdot) K_0 \pa R^u\\
   &+ \pa K_1 (\pz,\py) R^u + B
\end{split}
\end{equation}
where
\begin{equation*}
\begin{split}
           B               =& B_{II} + \pa y [\pt^2, S_1(z,y)]  R^u + B_{\textrm{Res}}\\
\ep^{11/2} B_{\textrm{Res}}=& (1 + Lz) \pt \ep^3 \patt + \pa y \pt^2 (K(z,y) \putt)\\
                            & +(1 + \pt^2 y) \pa \pt \py + \pt u \pa \putt.
\end{split}
\end{equation*}
The terms $B_{II}$ and $\pa y [\pt^2, S_1(z,y)]$
are small and smooth, and we can bound $B_{\textrm{Res}}$ {\it via} the
residual estimates, much as we did for $N^u_1$ above.  That is, we have
$$
\|B\|_{s-1} 
\le C\left( \ep^3  +  \ep^3 \|\bar{R}\|_{\sob^s_e} + \ep^{11/2} \|\bar{R}\|^2_{\sob^s_e} \right)
$$
under the same hypotheses as in the above Lemmas.  

At this time, it is tempting to simply invert 
$A(z,y)$.  Though we could do this, the
inverse of this operator is not smoothing.  In particular the presence of the term $\pa y K_0$
in $A$ will cause problems.  We can eliminate $K_0$ to highest order by letting 
$H_0(z,y) = (1 + Lz - \pa y K_0)$ act on (\ref{ishmael}).  We have for the first term
\begin{equation*}
\begin{split}
 &H_0(z,y) A(z,y) \pt R^a \\
=&(1+Lz)^2 \pt R^u + (1+Lz)\pa y K_0 \pt R^a-\pa y K_0\left( (1+L z) \pt R^a \right)\\
 &- \pa y K_0\left(\pa y K_0 \pt R^a\right) + H_0(z,y)(\pa y S_1(z,y) \pt R^a)\\
=&\left( (1+Lz)^2  + (\pa y)^2  + H_1(z,y)\cdot\right) \pt R^a
\end{split}
\end{equation*}
where
\begin{equation*}
\begin{split}
H_1(z,y)\cdot =& \pa y \left([Lz,K_0]\cdot - K_0 [\pa y,K_0] \cdot - (1+K_0^2) \pa y \cdot\right)\\ 
                &+ H_0(z,y)\left(\pa y S_1(z,y) \cdot \right).
\end{split}
\end{equation*}
Notice that $H_1$ is made up of smoothing operators, and is thus a smoothing operator.

Now, for the second term in (\ref{ishmael}) we have
\begin{equation*}
\begin{split}
 &H_0(z,y) (1 + \pt^2 y - \pt u K_0 \cdot) K_0 \pa R^u\\
=&(1 + Lz)(1 + \pt^2 y) K_0 \pa R^u - (1+Lz) \pt u K_0^2 \pa R^u\\ 
 &- \pa y K_0\left((1+\pt^2y) K_0 \pa R^u\right)+\pa y K_0\left(\pt u K_0^2 \pa R^u\right)\\
=&(1 + Lz)(1+\pt^2 y) K_0 \pa R^u\\
 & - \pa y [K_0,\pt^2y] K_0 \pa R^u + \pa y K_0\left(\pt u K_0^2 \pa R^u \right)\\
 &+\left( \pt u (1 + Lz) + \pa y (1 + \pt^2 y) \right) K_0^2 \pa R^u.
\end{split}
\end{equation*}
Notice that by comparing the last line of the above with
(\ref{ww 1 1}), we see that it is identically zero!  
One more rearrangement of this yields
\begin{equation*}
\left( (1 + Lz)(1 + \pt^2 y) - \pa y \pt u\right) K_0 \pa R^u + B_2
\end{equation*}
where
$$
B_2 = \pa y \left(-[K_0,\pt^2 y] + [K_0,\pt u] K_0 + \pt u (1+K_0^2)\right) K_0 \pa u
$$  
is a smooth and small function by Lemma \ref{commutator} in Section \ref{operator K 2}.  

If we let 
\begin{equation*}
\begin{split}
f =& ((1 + Lz)^2 + (\pa y)^2)^{-1} \\
g =& (1+Lz)(1+\pt^2 y) -\pa y \pt u
\end{split}
\end{equation*}
then we have transformed (\ref{ishmael}) into
\begin{equation*}
\begin{split}
0=&(1 + f H_1(z,y)) \pt R^a + f g K_0 \pa R^u\\
  &+ f H_0(z,y)\left(\pa\left(K_1(\pz,\py) R^u\right)\right)\\
  &+ f \left( B_2 +  H_0(z,y) B\right).
\end{split}
\end{equation*}

By the Neumann series,
$$
(1+f H_1(z,y)\cdot)^{-1} = 1 + \sum_{n=0}^{\infty} (-1)^n f^n H_1^n(z,y)\cdot.
$$
Since $H_1$ is smoothing, this is the identity plus a smoothing piece.  Let
$$
H_2(z,y)\cdot = \sum_{n=0}^{\infty} (-1)^n f^n H_1^n(z,y)\cdot.
$$
Thus,
\begin{equation*}
\begin{split}
0 = & \pt R^a + fg K_0 \pa R^u \\ 
    & + f H_0(z,y) \left(\pa\left(K_1(\pz,\py) R^u\right)\right) - N^a_1 - N^a_2
\end{split}
\end{equation*}
where
\begin{equation*}
\begin{split}
-N^a_1 = &H_2(z,y)
    \left(fg K_0 \pa R^u + f H_0(z,y) \left(\pa\left(K_1(\pz,\py) R^u\right)\right)\right) \\
-N^a_2 = &(1 + H_2(z,y)) f \left( B_2 +  H_0(z,y) B\right).
\end{split}
\end{equation*}

Finally we rewrite the above as
$$
\pt R^a = -\left(1 + K_0 a - L z + N^a_s\right)\pa K_0 R^u - \pa\left(K_1(\pz,\py) R^u\right) + N^a
$$
with
\begin{equation*}
\begin{split}
N^a_s = &f g - (1 + K_0 a - L z)\\
N^a = &N^a_1 + N^a_2 + N^a_3\\
N^a_4 =& -(f H_0(z,y) - 1) \left(\pa\left(K_1(\pz,\py) R^u\right)\right).
\end{split}
\end{equation*}
Notice that $(1+K_0 a- Lz)$ is the first order approximation to $f g$.  And thus,
using techniques exactly like those we used in proving the bounds on $N^y$ and $N^u$, we have
the following.
\begin{lemma}
For all $C_R>0$, there exists $\epsilon_0$ such that for all $\ep \in (0,\ep_0)$ 
and $t$ such that $\sup_{0 \le t^\prime \le t} \|\bar{R}_e(\cdot,t^\prime)\|_{\sob^s_e} \le C_R$
we have
$$
\max \left\{ \|N^a_s\|_{s-1},\|N^a\|_{s-1} \right\}
\le C\left( \ep^3  +  \ep^3 \|\bar{R}_e\|_{\sob^s_e} + \ep^{11/2} \|\bar{R}_e\|^2_{\sob^s_e} \right).
$$
\end{lemma}

Recapping, we have shown that the three dimensional system may be rewritten as
\begin{equation}
\label{three d error equation}
\begin{split}
\pt R^z &= K_0 R^u\\
\pt R^y &= K_0 R^u + M_1(\pz) \pa R^u  - (\py + K_0(\ep^2 \Psi^y K_0)) \pa R^u + N^y\\
\pt R^u &= -(1 - \ep^2 W_1) \pa R^y + N^u
\end{split}
\end{equation}
and the four dimensional system as
\begin{equation}
\label{four d error equation}
\begin{split}
\pt R^z =& K_0 R^u\\
\pt R^y =& K_0 R^u + M_1(\pz) \pa R^u \\& - (\py + K_0(\ep^2 \Psi^y K_0)) \pa R^u + N^y\\
\pt R^u =& R^a\\
\pt R^a =&-\left(1 + K_0 a - L z + N^a_s\right)\pa K_0 R^u\\& - \pa\left(K_1(\pz,\py) R^u\right) 
        + N^a.
\end{split}
\end{equation}

We remark now that these are only cosmetically different than the equations which 
determine the evolution of the error for the KdV approximation alone in 
\cite{schneider.etal:00}.   See {\it p} 1524 for the equations in three dimensions and
{\it p} 1526 in four.  Their variables 
$$
(Z_1,X_2,U_1,V_1)
$$ 
correspond to our 
$$
(z,y,u,a),
$$
and their functions 
$$
(N^2,N^3,N^4,N^5,N^8)
$$ 
are our 
$$
(N^z,N^y,N^u,N^a_s,N^a).
$$  
The only
difference of note is that their estimates contain a term they call $q(t)$ while ours do not.
This term, which is related  to the interaction of the left and right moving wavetrains,
has been removed in this paper by the inclusion of the function $W_3$ in the approximating
functions $\Psi$.  This simplification does not adversely affect
the means which they employ to prove that the error functions remain $O(1)$ over the long
time scale.   Therefore we appeal to their results on {\it pp} 1524-1533.  That is,
\begin{proposition}
\label{easy way out}
For all $T_0 >0 $, $s>4$ and $C_I>0$, there exists $\ep_0$ such that for all $0 \le \ep \le \ep_0$,
the unique solution $\bar{R}_e$ of (\ref{four d error equation}) 
with initial conditions such that 
$$
\|\bar{R}_e(\cdot,0)\|_{\sob^s_e}  \le C_I 
$$
satisfies
$$
\sup_{t \in [0,T_0 \ep^{-3}]} \|\bar{R}_e(\cdot,t)\|_{\sob^s_e} \le C
$$
where $C$ is independent of $\epsilon$.  
\end{proposition}

Implicit in the above Proposition
is the assumption that the initial conditions for the water wave problem have the form:
$$
\left( \begin{array}{c}
 z(\al,0)\\
 y(\al,0)\\
 u(\al,0)
\end{array}\right) 
=
\left( \begin{array}{c}
\pz(\al,0)\\ 
\py(\al,0)\\ 
\putt(\al,0) 
\end{array}\right) 
+ \ep^{11/2}\bar{R}_0(\al).
$$
So we see that this Proposition immediately proves Theorem 
\ref{main result water wave 2}.  

Now that we have this result, there are a few small steps,
and one big step, needed to prove Theorem \ref{main result water wave}.
The first simple step is to note that the $z$ is not a very physical coordinate 
and that we would  prefer estimates for $x_\al$.
Since $L$ is a bounded operator and gives the relationship
between both $z$ and $x_\al$ and $\pz$ and $\pd$, 
we have automatically 
$$
\sup_{t \in [0,T_0 \ep^{-3}]} \|x_\al(\cdot,t) - \pd(\cdot,t)\|_{H^s} \le C \ep^{11/2}.
$$

Secondly, the expressions for $\pd$, $\py$ and $\pu$ contain terms of $O(\ep^6)$.  These
terms were needed to make the residual sufficiently small, but they are unneccessary now.
Moreover, the appearance of the operator $L^{-1}$ and inverse derivatives
in the definitions of $\py$ and $\pu$ is not very intuitive.  And so, it is a simple
consequence of Lemma \ref{workhorse} and the triangle inequality that
\begin{equation*}
\begin{split}
\|\pd - \psi^d\|_s \le& C\ep^{11/2}\\
\|\py - \psi^y\|_s \le& C\ep^{11/2}\\
\|\pu - \psi^u\|_s \le& C\ep^{11/2}
\end{split}
\end{equation*}
where $\psi^d$, $\psi^y$ and $\psi^u$ were given in the Introduction in equations
(\ref{psi d 1}), (\ref{psi y 1}) and (\ref{psi u 1}). And so we have the corollary
\begin{corollary}
\label{XK}
If the initial conditions for (\ref{water wave}) are of the form
\begin{equation}
\label{useful i.c.}
\left( \begin{array}{c}
 x_\al(\al,0)\\
 y(\al,0)\\
 u(\al,0)
\end{array}\right) 
=
\left( \begin{array}{c}
\psi^d(\al,0)\\ 
\psi^y(\al,0)\\ 
\psi^u(\al,0) 
\end{array}\right) 
+ \ep^{11/2}\bar{R}_1(\al)
\end{equation}
with $\|\bar{R}_1\|_{\sob^s} \le C_I$
then the solution of (\ref{water wave}) satisfies the estimate:
$$
\left \|
\left( \begin{array}{c}
 x_\al(\cdot,t)\\
 y(\cdot,t)\\
 u(\cdot,t)
\end{array}\right) 
-
\left( \begin{array}{c}
 \psi^d (\cdot,t)\\
 \psi^y (\cdot,t)\\
 \psi^u (\cdot,t)
\end{array}\right) 
\right \|_{\sob^s}
\le C_F \epsilon^{11/2}
$$
for $t \in [0,T_0 \epsilon^{-3}]$.  The constant $C_F$ does not depend on $\epsilon$.
\end{corollary}

Finally, we must deal with initial conditions.  Recall from the discussion in
Section \ref{preliminaries} that it is typical to specify the initial data
for the water wave problem in the long wavelength, small amplitude limit
by
\begin{equation}
\label{typical i.c.}
(\bar{x}_{\ab}(\ab,0),\bar{y}(\ab,0),\bar{u}(\ab,0)) = (0,\ep^2 \Theta_y(\ep \ab),\ep^2 \Theta_u(\ep \ab)).
\end{equation}
However, the above results are applicable if the initial data is of the form
seen in (\ref{useful i.c.}).
We eliminate this discrepency by altering the initial parameterization of the free surface.  
What should this change be?  Clearly,
\begin{equation}
\label{hope for}
\begin{split}
\ab= \al + x(\al,0).
\end{split}
\end{equation}

Now set 
$U(\beta,0) = U_0(\beta)$,
$V(\beta,0) = V_0(\beta)$,
$F(\beta,0) = F_0(\beta)$,
$G(\beta,0) = G_0(\beta)$,
and
$P(\beta,0) = 0$, and let
\begin{align*}
\ab = &\al + \int_0^{\al} \psi^d(a,0) da\\
    = &\al + \ep X_1(\ep \al) + \ep^3 X_2 (\ep \al)
\end{align*}
where
\begin{equation}
\label{X_1 X_2}
\begin{split}
\ep   X_1(\al) &= -\ep^2 \int_0^{\al} \left(U_0(\ep a) + V_0(\ep a)\right) d a\\
\ep^3 X_2(\al) &= -\ep^4 \int_0^{\al} \left(F_0(\ep a) + G_0(\ep a)\right) d a.
\end{split}
\end{equation}
With this definition, we clearly have satisfied the first condition in
(\ref{useful i.c.}).
We also want
\begin{equation*}
\begin{split}
\Theta_y(\ep \ab) = \ep^2 \psi^y(\al) + O(\ep^{11/2})\\
\Theta_u(\ep \ab) = \ep^2 \psi^u(\al) + O(\ep^{11/2})
\end{split}
\end{equation*}
or rather
\begin{equation*}
\begin{split}
\ep^2\Theta_y(\ep \al + \ep^2 X_1(\ep \al) + \ep^4 X_2(\ep \al)) = \ep^2 \psi^y(\al) + O(\ep^{11/2})\\
\ep^2\Theta_u(\ep \al + \ep^2 X_1(\ep \al) + \ep^4 X_2(\ep \al)) = \ep^2 \psi^u(\al) + O(\ep^{11/2}).
\end{split}
\end{equation*}
Applying Taylor's theorem we have
\begin{equation*}
\begin{split}
 \Theta_y + \ep^2 X_1\Theta_y^{\prime}
=&(U_0 + V_0) + \ep^2 (\third \pbm^2 U_0 + \third \pbp^2 V_0)\\& + \ep^2(F_0 + G_0) + \ep^2 (U_0+V_0)^2 \\
\Theta_u + \ep^2 X_1\Theta_u^{\prime}
=&(U_0 - V_0) + \ep^2 (\frac{1}{6} \pbm^2 U_0
 - \frac{1}{6} \pbp^2 V_0)\\& + \ep^2(F_0 - G_0) + \ep^2 (\frac{3}{4}U_0^2-\frac{3}{4}V_0^2) .
\end{split}
\end{equation*}

We can solve the above by taking
\begin{equation}
\label{kdv i.c.}
\begin{split}
U_0 = 1/2 (\Theta_y + \Theta_u)\\
V_0 = 1/2 (\Theta_y - \Theta_u)
\end{split}
\end{equation}
and
\begin{equation}
\label{lkdv i.c.}
\begin{split}
F_0 = 1/2 (h_y + h_u)\\
G_0 = 1/2 (h_y - h_u)
\end{split}
\end{equation}
where
\begin{align*}
h_y &= X_1 \Theta_y^\prime -  \third \pbm^2 U_0 - \third \pbp^2 V_0 - (U_0+V_0)^2\\
h_u &= X_1 \Theta_u^\prime - \frac{1}{6} \pbm^2 U_0
      + \frac{1}{6} \pbp^2 V_0 -\frac{3}{4} U_0^2+\frac{3}{4}V_0^2.
\end{align*}

The functions $U_0$, $V_0$, $F_0$ and $G_0$ are all in $H^s(4)$,
and so the use of Taylor's theorem is justified by Lemma \ref{mean value theorem}.  
So we have proven:
\begin{lemma}
Given initial conditions for the water wave equation in the form (\ref{typical i.c.}),
define $U_0$, $V_0$, $F_0$, $G_0$, $X_1$ and $X_2$ as in (\ref{kdv i.c.}), (\ref{lkdv i.c.})
and (\ref{X_1 X_2}).   Then the 
reparameterization of the initial profile given by:
$$
\ab = \al + \ep X_1(\ep \al) + \ep X_2(\ep \al)
$$
results in initial conditions given by (\ref{useful i.c.}).
\end{lemma}

\begin{remark}
Let $\varphi^\pm(\beta_\pm,0) = \varphi^\pm_0(\beta_\pm)$.  Then this
Lemma will still be true if we replace $F_0$ with $\varphi^-_0$ and 
$G_0$ with $\varphi^+$ and set $F_0$ and $G_0$ to be identically zero.
That is, we have some choice in the way we select the initial conditions for
the higher order equations.  
\end{remark}

Combining this Lemma with Corollary \ref{XK}
we prove Theorem \ref{main result water wave}.  So we are done.

\section{Assorted proofs}
\label{proofs}

\begin{proof} {\bf For Lemma \ref{inverse lemma}}:
Let $f^{-1}(\xi) = \xi - g_2(\xi)$.  Since $f^{-1}(f(\al)) = \al$ we have
$$
\al = f(\al) - g_2(f(\al)),
$$
or rather
\begin{equation}
\label{master inverse}
g_2(f(\al)) =  g(\al)
\end{equation}
Notice that this relation implies  $\|g_2\|_{L^{\infty}} = \|g\|_{L^\infty}$.  Taking
a derivative, we have
$$
g_2^\prime(f(\al)) = \frac{g^\prime(\al)}{1 + g^\prime(\al)}
$$
which implies that 
$\|g_2^\prime\|_{L^{\infty}} \le C \|g^\prime\|_{L^\infty}$.  
If we expand the left hand side of (\ref{master inverse}) by the mean value theorem
we see
\begin{align*}
g_2(\al + g(\al)) &=  g(\al)\\
g_2(\al) + g_2^\prime (\al^*)g(\al) &=  g(\al).
\end{align*}
This implies
$g_2(\al) = g(\al) + O(\|g^\prime\|_{L^\infty} \|g\|_{L^\infty})$.
Now, (\ref{master inverse}) can be rewritten and expanded using Taylor's theorem:
\begin{align*}
g_2(\xi) &= g(f^{-1}(\xi))\\
         &= g(\xi - g(\xi) + O(\|g^\prime\|_{L^\infty} \|g\|_{L^\infty}))\\
         &= g(\xi) + g^{\prime}(\xi) (- g(\xi)+ O(\|g^\prime\|_{L^\infty} \|g\|_{L^\infty}))\\
         &+ 1/2 g^{\prime \prime}(\xi^*)(- g(\xi)+ O(\|g^\prime\|_{L^\infty} \|g\|_{L^\infty}))^2\\
         &= g(\xi) - g(\xi) g^{\prime}(\xi) - E
\end{align*}
which completes the proof.  
\end{proof}

\begin{proof} {\bf For Lemma \ref{the weight}}:
Since $(1 + \beta^2)^{n/2} F(\beta) \in H^s$, by 
the Sobolev embedding Theorem there
is a $C$ such that
$$
F(\beta) \le C(1 + \beta^2)^{-n/2}.
$$
So
\begin{equation*}
\begin{split}
     &\int_{|\al| > C_0 \ep^{-3}} |F(\ep \al)|^2 d\al\\
\le &C\int_{|\al| > C_0 \ep^{-3}} |1 + (\ep \al)^2|^{-n} d\al\\ 
\le &C\int_{|\al| > C_0 \ep^{-3}} |\ep \al|^{-2n} d\al\\ 
\le &C\ep^{-2n} \int_{|\al| > C_0 \ep^{-3}} |\al|^{-2n} d\al\\ 
\le &C\ep^{-2n} (\ep^{-3})^{-2n+1}\\
\le &C\ep^{4n-3}.
\end{split}
\end{equation*}
The higher derivatives are bounded in exactly the same fashion.  The extra powers of
$\ep$ come from the long wavelength scaling.
\end{proof}

\begin{proof} {\bf For Lemma \ref{workhorse}}:
The proof is a straightforward calculation.  
\begin{eqnarray*}
&   & \|A f(\cdot) - A_n f(\cdot)\|^2_s\\
&=  & \int (1+k^2)^s  |(\widehat{A}(k)-\widehat{A_n}(k))\widehat{f}(k)|^2 dk\\
&\le& C\int (1+k^2)^s  |k^n \widehat{f}(k)|^2 dk\\
&=  & C\int (1+k^2)^s  |\widehat{\partial^n_x f}(k)|^2 dk\\
&=  & C\|\partial^n_x f(\cdot)\|^2_s.
\end{eqnarray*}
The proof for long wavelength data follows immediately from this.
\end{proof}

\begin{proof} {\bf For Lemma \ref{geometric}}:
The fact that $g(\ep x)$ is bounded as such in $L^2$ follows automatically from
the geometric series approximation.  That is, since $|\ep^2 f(\ep x)| \le 1/2$,
we know that:
$$
|g(\ep x)| \le C|\ep^2 f(\ep x)|^2.
$$
And thus we have:
\begin{align*}
&\|g(\ep \cdot)\|_{L^2}^2 \\
&\le C\int \ep^4 |f(\ep x)|^4 dx\\
&\le C\ep^4 \|f(\cdot)\|^2_{L^\infty} \int f^2(\ep x) dx\\
&\le C\ep^{7/2} \|f(\cdot)\|_s.
\end{align*}
Now consider the $L^2$ norm of $g^\prime(\ep x)$.  A direct calculation shows that
$$
\frac{d}{dx} g(\ep x) = - \ep^3 f^\prime(\ep x) 
                             \left(\frac{1}{(1+\ep^2 f(\ep x))^2} -1 \right).
$$
Taylor's theorem shows that 
$$
|\frac{d}{dx} g(\ep x)| \le C|\ep^3 f^\prime(\ep x)||\ep^2 f(\ep x)|.
$$
And so, just as before we have that this is bounded by $C \ep^{9/2}$ (which is of course
bounded by $C \ep^{7/2}$).  

We could keep on going in this fashion---showing each derivative of $g$ is bounded.  This
is however difficult as finding higher and higher derivatives is a notationally taxing
job---see the expression of Faa-di-Bruno  for proof of that!  

Instead we take the following
approach.  Let
$$
h(y) = \frac{1}{1+y} - 1 + y.
$$
For $y \in [-1/2,1/2]$, this function is real analytic and there exists another function
$\tilde{h}(y)$ (real and analytic on the same interval) such that $h(y) = y \tilde{h}(y)$.
Now, define $\tilde{h}_\ep(Y) = \tilde{h}(\ep^2 y)$.  We have that:
$$
\|\tilde{h}_\ep(\cdot)\|_{C^s[-1/2 \ep^{-2}, 1/2 \ep^{-2}]} 
\le  
\|\tilde{h}(\cdot)\|_{C^s[-1/2 , 1/2 ]}.
$$ 
The point here is that the $C^s$ norm of $\tilde{h}_\ep$ can be bounded 
{\it independently} of $\ep$.  

Now notice that $g(X) = \ep^2 f(X) \tilde{h}_\ep(f(X))$.  Since 
$f \in H^{s+1}$, we know that $f \in C^s$.  This implies that 
$\tilde{h}_\ep(f(X)) \in C^s$, with $C^s$ norm bounded independent of
$\ep$.  Thus we have $f(X) \tilde{h}_\ep(f(X)) \in H^s$, with
a norm bounded independent of $\ep$.  With this in hand,
we have that $\|g(\cdot)\|_s \le C \ep^2$, with $C \ne C(\ep)$.  Now,
the derivatives of $g$ can be bounded as follows:
\begin{align*}
&\| \frac{d^n}{dx^n} g(\ep \cdot) \|_{L^2}\\ 
=&\ep^n \|g^{(n)} (\ep \cdot) \|_{L^2}\\
\le &C \ep^{n-1/2} \|g(\cdot)\|_s\\
\le &C \ep^{n+3/2}.
\end{align*}
Provided $n \ge 2$, this term is small enough.  And so we have shown that
the first $s$ derivatives are sufficiently small in $L^2$, and we have proved the
Lemma.
\end{proof}

\nocite{*}
\bibliographystyle{plain}
\bibliography{ref_thesis_old.bib}

\end{document}